\documentclass[10pt]{article} 
\usepackage{amsfonts}  
\usepackage{amssymb}
\usepackage{mathrsfs}
\usepackage{latexsym}

\title{Constructive Coordinatization of Desarguesian Planes}

\author{Mark Mandelkern}

\begin{document}

\newtheorem{defn}{Definition}[section]   
\newtheorem{prop}[defn]{Proposition}
\newtheorem{lm}[defn]{Lemma}
\newtheorem{thm}[defn]{Theorem}
\newtheorem{cor}[defn]{Corollary}
\newtheorem{prob}[defn]{Problem}
\newtheorem{exer}[defn]{Exercise}
\newtheorem{ex}[defn]{Example}
 
\setlength{\parindent}{15pt}

\begin{center}
\begin{Large}
Constructive Coordinatization of Desarguesian Planes
\end{Large}\\

\vspace*{5mm}
Mark Mandelkern \\ 
\end{center}
\vspace*{1mm}

\pagenumbering{arabic}
\setcounter{page}{1}

\begin{quotation}
 \noindent \textbf{Abstract.} A classical theory of Desarguesian geometry, originating with D. Hilbert in his 1899 treatise \textit{Grundlagen der Geometrie,} leads from axioms to the construction of a division ring from which coordinates may be assigned to points, and equations to lines; this theory is highly nonconstructive. The present paper develops this coordinatization theory constructively, in accordance with the principles introduced by Errett Bishop in his 1967 book, \textit{Foundations of Constructive Analysis}.

The traditional geometric axioms are adopted, together with two supplementary axioms which are constructively stronger versions of portions of the usual axioms. Stronger definitions, with enhanced constructive meaning, are also selected; these are based on a single primitive notion, and are classically equivalent to the traditional definitions. Brouwerian counterexamples are included; these point out specific nonconstructivities in the classical theory, and the consequent need for strengthened definitions and results in a constructive theory.

All the major results of the classical theory are established, in their original form --- revealing their hidden constructive content. 
\end{quotation}

\hspace*{5mm} \\

\section{Introduction}
\label{sec:intro}

In various forms, the constructivist program goes back to Leopold Kronecker, Henri Poincar\'e, L. E. J. Brouwer, and many others. The most significant recent work, using the strictest methods, is due to Errett Bishop. A large portion of analysis has been constructivized in Bishop's book, \textit{Foundations of constructive analysis} [B]; this book also serves as a guide  for constructive work in other fields.\footnote{Expositions of constructivist ideas and methods, and further references, may be found in [B, BB; Chapter 1], [BM], [R], [S], and [MRR; pp. 1-6].} 

The initial phase of this program involves the rebuilding of classical theories using only constructive methods. This phase is based on the entire body of classical mathematics, as a wellspring of theories waiting to be constructivized. ``Every theorem proved with [nonconstructive] methods presents a challenge: to find a constructive version, and to give it a constructive proof.''\footnote{Errett Bishop, 1967 [B, page x].} 

Following this dictum, the present work is based on the classical theory of Desarguesian planes and their coordinatization, which originated with D. Hilbert [H]; the plan used here follows the modern presentation given by E. Artin [A]. The classical theory is highly nonconstructive; it relies heavily, at nearly every juncture, on the \textit{principle of the excluded middle}. For example, it is assumed that a given point is either on a given line, or not on the line --- although no finite routine is provided for making such a determination. 

The Desarguesian coordinatization theory will be developed constructively, adhering to the precepts put forth by Bishop. The selection of axioms and definitions is constrained by the constructive properties of the real plane $\mathbb{R}^{2}$; we will not expect to prove any theorem that is constructively invalid on $ \mathbb{R}^{2} $. 

 \textit{Primitive notion.} We adopt a  single primitive notion, ``distinct points,'' with strong properties. The most significant  property involves a disjunction. Classically, the condition is obvious: ``A point cannot be equal to each of two given distinct points.'' Constructively, we must have a finite routine that results in a definite decision: The point is distinct from one of the two given points --- which point? 

 \textit{Principal relation.} Rather than the usual  concept ``point on a line,''  the concept ``point outside a line'' will have constructive utility as the principal relation. This relation is given an affirmative definition directly in terms of the primitive. 

 \textit{Axioms.} The axioms fall into three groups. The five classical axioms used by Artin [A] are adopted here without change. In group \textbf{G} are the three traditional axioms for plane geometry. In group \textbf{K} are the two symmetry axioms required for the coordinatization.  

 \textit{Axiom group \textbf{L.}} The two axioms in this group are not essentially new, but are inherent in the classical theory. These axioms are strong versions of the uniqueness portions of the axioms in group \textbf{G}; they provide disjunctive decisions in situations involving nonparallel lines. One of these supplementary axioms follows, in classical form, from the parallel postulate: ``A line cannot be parallel to each of two given, distinct, intersecting lines.''  Constructively, a finite routine must provide a definite decision: The line is nonparallel to one of the given lines --- which line? The alternative forms of these axioms will be discussed further in appendix A to section \ref{sec:ax}.

 \textit{Parallelism.} This is the central idea of the theory. The  affirmative concept, ``nonparallel,'' is the focus; in turn, this concept will be based on an affirmative definition of ``distinct'' lines.

\textit{Dilatations.} The coordinatization is based on the symmetries of the geometry; these are the dilatations, maps that  preserve direction. The classical theory of dilatations rests heavily on nonconstructive principles. Here we must strengthen the definition, requiring a dilatation to be injective in a strict sense. New constructions are required to show that the inverse of a dilatation is also a dilatation, and to prove an extension theorem for dilatations. 

 \textit{Translations.} The classical definition of a translation, a dilatation that is either the identity or has no fixed point, is not constructively feasible. We will say that a dilatation is a translation if any traces are parallel.\footnote{This is classically equivalent to the traditional definition, but only under the slightly-modified definition of trace that is used here.} Then it must be proved that a translation that maps one point to a distinct  point has the same behaviour at every point.
 
\textit{Coordinatization.} The scalars used for the coordinatization are certain homomorphisms of the translation group. A new construction is required to show that non-zero homomorphisms are injective, so that the scalars form a division ring. 

 \textit{Desargues's Theorem.} This theorem will be shown equivalent to the symmetry axioms. Desargues's Theorem may thus be used as an alternative to these axioms; it has the advantage that it involves only direct properties of the parallelism concept. The proof requires new constructions, and the extension theorem for dilatations. 

 \textit{Pappus's Theorem.} It will be shown that commutativity of the division ring  is equivalent to Pappus's Theorem. The proof is based on the preceding work, and requires no new constructions.

 \textit{Geometry based on a field.} From a  given field with suitable properties, we  construct a Desarguesian plane. The classical theory does this for an arbitrary division ring; constructively, the case of a  (non-commutative) division ring is an open question. 

 \textit{The real plane.} The constructive properties of the field $ \mathbb{R} $ of real numbers  ensure that $ \mathbb{R}^{2} $ is a Desarguesian plane satisfying all the axioms. The order and metric structures on the reals allow possible alternatives for the principal relation, ``point outside  a line.'' These alternatives are shown to be equivalent to the adopted definition. 

 \textit{Brouwerian counterexamples.} These counterexamples pinpoint the nonconstructivities  of the classical theory, and facilitate the selection of axioms and definitions for the constructive theory. Two of the examples show that the following statements are constructively invalid: ``Parallel lines are either equal or disjoint.'' ``If two lines have a unique point in common, then they are nonparallel.''
 
 \textit{Other constructive geometries.} Constructive geometry has been approached from various directions. The work of A. Heyting [He1] concerns projective geometry, obtaining a coordinatization by means of projective collineations, whereas the present paper is a constructive study of parallelism. The paper [He2] concerns axioms for plane incidence geometry and extensions to projective planes.\footnote{A few comments on Heyting's axiom system will be given in appendix B to section \ref{sec:ax}.} Interesting papers by D. van Dalen [D1, D2] concern primitive notions and relations between projective and affine geometry.

Other work is more closely related to logic, type theory, recursive function theory, and computer techniques --- approaches far removed from the``straightforward realistic approach" [B, p.10] proposed by E. Bishop and followed in the present paper. The interesting papers by J. von Plato [P1, P2], D. Li, X. Li, P. Jia [L, LJL], Lombard and Vesley [LV] enable valuable comparisons between the several varieties of constructivism (see [BR]), in the context of geometry. 

\textit{Logical setting.} This work uses informal intuitionistic logic; it does not operate within a formal logical system.\footnote{For presentations of informal intuitionistic logic, as it is used in Bishop-type constructive mathematics, see [BR] and [BV]. For the most part, it suffices that one exercise assiduous restraint in regard to the connective ``or''.} This constructivist principle has been most concisely expressed as follows: ``Constructive mathematics is not based on a prior notion of logic; rather, our interpretations of the logical connectives and quantifiers grow out of our mathematical intuition and experience.''\footnote{[BR, page 11]} For the origins of modern constructivism, and the disengagement of mathematics from formal logic, see Bishop's Chapter 1, ``A Constructivist Manifesto'', in [B].

 \textit{Summary.} The entire coordinatization theory of Desarguesian planes has been constructivized in the spirit of Bishop-type constructivism. The results retain their original classical form, with enhanced constructive meaning. A number of open problems remain; these have been noted in the various sections. 

\section{Axioms}
\label{sec:ax}

A geometry will at first consist of a set of points, a set of lines, and a single primitive relation. The principal concepts will be defined in terms of this primitive. The first two axiom groups will be introduced: the traditional group \textbf{G}, and the special group \textbf{L} concerning nonparallel lines. The elementary properties of parallelism will be derived.
 
\begin{defn} 
\label{ax2}
\textnormal{A} geometry 
\textnormal{ $ \mathscr{G} = ( \mathscr{P}, \mathscr{L} ) $ consists of:\\
\hspace*{5mm}(a) A set $\mathscr{P}$, whose elements are called  ``points,'' with a given equality relation.\\
\hspace*{5mm}(b) A set $\mathscr{L}$ of subsets of $ \mathscr{P} $, called ``lines,'' with the usual equality relation for subsets. When  $ P \in l $, we say that \textit{``The point $ P $ lies on the line $ l $,''}  
or that \textit{``The line $ l $ passes through the point $ P $.''}\\
\hspace*{5mm}(c) An inequality 
relation on the set of points $ \mathscr{P} $, 
written $ P \ne Q $; we say that \textit{``The points $ P $ and $ Q $ are distinct.''}  This relation is invariant with respect to the equality relation on $\mathscr{P}$, and has the following properties:\\
\hspace*{10mm} (c1) $  \neg (P \ne P) $\\
\hspace*{10mm} (c2) If $  P \ne Q $, 
then $  Q \ne P $.\\ 
\hspace*{10mm} (c3) If $ P $ and $ Q $ are distinct points, then any point $ R $ is either distinct 
from $ P $ or distinct from $ Q $.\\
\hspace*{10mm} (c4) If $ \neg (P \ne Q ) $, 
then $ P = Q $.}
\end{defn}

The converse of condition \ref{ax2}(c4) follows from  condition (c1). However, the statement ``If $ \neg (P = Q ) $, then $ P \ne  Q $'' is constructively invalid on the real plane $ \mathbb{R}^{2} $.\footnote{See example \ref{brouA}.} Thus the notation $ P \ne Q $ is not used in the usual classical sense of a negation. The principal relation, $ P \notin l $, will also have an affirmative meaning: 

\begin{defn}
\label{ax-2b}
\textnormal{We define a relation between the points of $ \mathscr{P} $ and the lines of $ \mathscr{L} $ as follows: 
\begin{displaymath}
P \notin l  \textnormal{~~~if~~~}  P \ne Q \textnormal{~~for all points~} Q \in  l   
\end{displaymath}
We say that \textit{``The point $ P $ lies outside the line $ l $.''}}
\end{defn} 

\begin{prop}
\label{outside1} \hspace*{5mm}\\ 
\hspace*{5mm}(a) The relation ``$ P \notin l $'' is invariant with respect to the equality relations on $ \mathscr{P} $ and $ \mathscr{L} $.\\
\hspace*{5mm}(b) If $ P \in l$, then $ \neg (P \notin l) $.

\end{prop}

The converse of condition \ref{outside1}(b) will be most essential; it will be established in theorem \ref{outside2}, after the first two axiom groups are introduced.\\ 

\noindent \textit{Notes.} It is traditional to define a geometry so that the lines are independent of the set of points, rather than as subsets of points. It is not difficult to do this constructively. There results, as usual, a correspondence between lines and sets of points; an equivalent geometry may be formed in which the lines are, in fact, sets of points. Thus it is expedient to simply define lines as subsets of the given set of points, as in definition \ref{ax2}. 

 The properties of the primitive relation $ P \ne Q $ and the principal relation $ P \notin l $  have analogues in the constructive properties of  the real field $\mathbb{R}$ and  the real plane $\mathbb{R}^{2}$. On the line, the condition $|x|>0$ is affirmative; the condition $x=0$ is its negation. The statement ``$ \neg (x=0) $ implies $|x|>0 \,$'' is constructively invalid.\footnote{See section \ref{sec:brou}.} On the real plane $\mathbb{R}^{2}$, the basic relations correspond to the distance between points, and from a point to a line.\footnote{See section \ref{sec:real}.} Thus the statements ``$ \neg (P = Q) $ implies $ P \ne Q $'' and ``$ \neg (P \in l) $ implies $ P \notin l $'' are constructively invalid on the real plane.\footnote{See example \ref{brouA}.}

 Condition \ref{ax2}(c3) may be compared to the constructive dichotomy lemma for the real numbers: ``If $ a < b $, then for any $ x $, either $ x < b $ or $ x > a $.'' This lemma  serves for the constructive development of analysis in lieu of the classical trichotomy, which is constructively invalid [B, BB]. The validity of condition (c3) on the  real plane $\mathbb{R}^{2}$ results from applying the dichotomy lemma to coordinates of points. 

 The construction of a geometry according to definition \ref{ax2} must include algorithms, or finite routines, for the conditions listed; the same rule applies to the definitions and axioms that follow. The notion of algorithm is taken as primitive; for a discussion of finite routines and algorithms, see [BR, Chapter 1].\\

 \noindent \textbf{Notation and conventions.} 
\label{ax6} 
For any lines $ l $ and $ m $, the expression $ l \cap m \ne \emptyset $ will mean that there exists a point $ P $ such that $ P \in l \cap m $; i.e., there exists a finite routine that would produce the point $ P. $ When $ \mathscr{M} \subseteq \mathscr{L}$, the expression $ \mathscr{M} \ne \emptyset $  will mean that there exists a line $ l $  that is in the set $ \mathscr{M.}$ The expression ~$ = \emptyset $~ will mean  that the condition  ~$ \ne \emptyset $~ for the set in question leads to a contradiction. These conventions ensure that the expression ~$ \ne  \emptyset $~  does not mean merely that it is contradictory that the set in question is void. The symbol ~$ \equiv $~ will be used to define objects.

An inequality relation\label{apart} that satisfies conditions (c1) through (c3) of definition \ref{ax2} is called an \textit{apartness}; if it satisfies condition (c4) it is said to be \textit{tight}. For a comprehensive treatment of constructive inequality relations, see [MRR; \S\,I.2].

For maps \label{maps} between sets each having a tight apartness as an inequality relation, the usual equality and inequality relations will be used. Thus \label{sigma=} $ \varphi = \psi $ if $ \varphi x = \psi x $ for all $ x ,$ while $ \varphi  \ne \psi $ if there exists at least one $ x $ such that $ \varphi x \ne \psi x $. It follows that the inequality relation on a set of maps is also a tight apartness. 

A map  $ \varphi $ will be called\label{inj} \textit{injective} if $ x \ne y $ implies $ \varphi x \ne \varphi y $. The condition normally used in classical work, $ \varphi x = \varphi y $ implies $ x = y $, is called \textit{weakly injective}; although classically equivalent to ``injective,'' this condition is constructively far weaker, and is of minimal use here. A \textit{bijection} is injective and onto, and has an injective inverse.\\ 

\noindent \textbf{Parallelism.} \label{ax7}
The usual classical definition, two lines are parallel if they are either equal or disjoint, is constructively invalid on the real plane $ \mathbb{R}^{2} $.\footnote{See example \ref{brouB}.} From the various classically equivalent conditions, and conditions for ``nonparallel,'' we  select the strongest form of ``nonparallel'' as a definition, and then take ``parallel'' as the negation. In turn, the concept of ``nonparallel'' lines will depend on the concept of ``distinct'' lines, defined  below. Constructive difficulties arise if different, albeit classically equivalent, definitions are used.\footnote{See example \ref{brouH}.}

\begin{defn}
\label{ax4}
\textnormal{We define an inequality relation on
$ \mathscr{L}$ as follows: $ l \ne m $  if  there exists a point $ P \in l $ with $ P \notin m $, or if there exists a point $ Q \in m $ with $ Q \notin l $. We say that  \textit{``The lines $ l $ and $ m $ are distinct.''}}
\end{defn}

\begin{prop}
\label{lines1} \hspace*{5mm}\\
\hspace*{5mm}(a) The relation ``$ l \ne m $'' is  invariant with respect to the equality relation on $ \mathscr{L} $.\\
\hspace*{5mm}(b) $ \neg (l \ne l) $\\
\hspace*{5mm}(c) If $  l \ne m $, 
then $ m \ne l $.
\end{prop} 

Additional properties of this relation will be given in proposition \ref{lines2}, after the first two axiom groups are introduced.   

\begin{defn}
\label{ax8}
\textnormal {We define a relation $ \nparallel $ 
on $ \mathscr{L} $  as follows: 
\begin{displaymath}
  l \nparallel m \qquad \mathrm{if} \qquad 
 l \not=m  ~~\mathrm{and}~~  l \cap m \ne \emptyset 
\end{displaymath} 
We  say that \textit{``The lines l and m are nonparallel.''}\\
\hspace*{5mm}When $ \neg ( l \nparallel m )$, 
we write $ l \parallel m $, and say that 
\textit{``The lines l and m are parallel.''}}
\end{defn}

\begin{prop}
\label{par-inv}
The relations ``parallel'' and ``nonparallel'' are  invariant with respect to the equality relation on $ \mathscr{L} $. 
\end{prop}

It will be shown in proposition \ref{ax28} that the relation ``nonparallel'' is invariant with respect to the relation ``parallel.'' It will be shown in proposition \ref{ax29} that ``parallel'' is an equivalence relation.\\ 

 \noindent \textbf{Axiom groups.} The \label{ax9} axioms required for a \textit{constructive Desarguesian plane}\label{despl} $\mathscr{G}$ fall into three  groups. In axiom group \textbf{G} are the first three of the usual axioms for plane geometry. In group \textbf{L} are axioms concerning nonparallel lines. Group \textbf{K} will be introduced in 
sections \ref{sec:ring} and \ref{sec:coor}, to enable the coordinatization. The axioms in group \textbf{K} are equivalent to \textit{Desargues's Theorem}; this will be shown in section \ref{sec:des}. 

 The five axioms in groups \textbf{G} and \textbf{K} are virtually identical to those used classically in [A]. The axioms in group \textbf{L} are inherent in the classical theory; in this sense, no new axioms are introduced.\\

\noindent \textbf{Axiom group G.}  
\label{ax10} 
Although these axioms are the same as those used classically, their meanings are strengthened by the stronger definitions adopted here.\\

\noindent \textbf{Axiom G1.} \textit{Let $ P $ and $ Q $ be distinct points. Then there exists a unique line $ l $ such that the points $ P $ and $ Q $ both lie on $ l $.}

\begin{defn}
\label{ax11}  
\textnormal{We  denote the line generated in axiom G1 by $ P + Q $. Thus the statement $ l $ = $ P + Q $ will include the covert condition $ P \ne Q $.}
\end{defn}

\noindent \textbf{Axiom G2.}\label{parpost} \textit{Let $ P $ be any point and let $ l $ be any line. Then there exists a unique line $ m $ through $ P $ that is parallel to $ l $.}\\

\noindent \textbf{Axiom G3.} \textit{There exist three non-collinear points. That is, there exist distinct points $ A $, $ B $, $ C $ such that $ C \notin A+B $.}\\

\noindent \textbf{Axiom group L.} 
\label{ax14}  
These axioms are classical equivalents of the uniqueness portions of the axioms in group \textbf{G}; see the appendix to this section. To state axiom L1, we require first a proposition. 

\begin{prop}
\label{ax15}
Let $ l $ and $ m $ be nonparallel lines. Then there exists a unique point $ P $ such that $ P \in l \cap m $.
\end{prop}  

\noindent \noindent Proof. By definition \ref{ax8}, we have at least one point $ P $ in $ l \cap m $. Now 
let $ Q $ be any point in $ l \cap m $.  Suppose that 
$ Q \ne P $; it then follows from  axiom G1 that
$ l = m $, a contradiction. Thus $ \neg (Q \ne P) $, and by definition \ref{ax2} this means that $ Q = P $. $\Box$\\

The \label{ax16} converse: ``If ~$ l \cap m $~ consists of exactly one point, then $ l \nparallel m,$'' is constructively invalid on the real plane $\mathbb{R}^{2}$.\footnote{See example \ref{brouD}.} 

\begin{defn}
\label{ax18}
\textnormal{For the unique point of intersection  determined in proposition 
\ref{ax15}, we write simply $ P = l \cap m $.}
\end{defn}

\noindent \textbf{Axiom L1.} \textit{Let $ l $ and $ m $ be nonparallel lines, and let 
$ P $ be the point of intersection. Then for any 
point $ Q $ distinct from $ P $, 
either $ Q $ lies outside $ l $, or $ Q $ lies outside $ m $}.\\

\noindent \textbf{Axiom L2.} \textit{Let $ l $ and $ m $ be nonparallel lines. Then for any line $ n $, 
either $ n $ is nonparallel to $ l $, or $ n $ is nonparallel to $ m $.}\\

 \noindent \textit{Problem.} Find a single axiom to replace axioms L1 and L2.

\begin{prop}
\label{ax28}
The relation ``nonparallel'' is invariant with respect to the relation ``parallel.''
\end{prop}

\noindent Proof. Let $ l $ and $ m $ be nonparallel lines, and let $ n $ be a line parallel to $ m $. It follows from axiom L2 that either $ n $ is nonparallel to $ l $,  or $ n $ is nonparallel to $ m $; hence $ n $ is nonparallel to $ l $.\footnote{This conclusion follows because the other case contradicts the hypothesis. To some eyes, this sort of argument may appear as a ``proof by contradiction,'' contrary to a proper constructivist attitude. However, this method  merely involves the  ruling out of a case that does not occur. For further comment on this issue, see [B, Appendix B].} $\Box$

\begin{thm}
\label{outside2}
Let $ P $ be any point, and let $ l $ be any line. If $ \neg (P \notin l) $, then $ P \in l$.
\end{thm}

\noindent Proof. Let it be given that $ \neg (P \notin l) $. Axiom G3 provides a pair of nonparallel lines.  Using axiom L2, we find that one of these lines is nonparallel to $ l $; denote it by $ m $. Use the parallel postulate to construct the line $ n $ through $ P $ that is parallel to $ m $. It follows from proposition \ref{ax28} that $ n $ is also nonparallel to $ l $; set $ R \equiv l \cap n $. 

Suppose that $ P \ne R $. It then follows from axiom L1 that either $ P \notin l $, or $ P \notin n $. The first case is ruled out by our hypothesis; the second case is ruled out by the choice of $ n $. This contradiction shows that $ P = R $. Hence $ P \in l $. $\Box$

\begin{thm} 
\label{ax20}
If three non-collinear points are given as in axiom G3, then the three lines formed are nonparallel in pairs.
\end{thm}

\noindent Proof. Since $ A + B $ and $ A + C $ have the common point A, and $ C $ $ \notin $ $ A + B $, we have $ A + B \nparallel A + C $. 

 Since $ A + B  \; \cap \; A + C = A $, and $ B \ne A $, it follows from  axiom L1 that $ B \notin A + B $ or $ B \notin A + C $; thus $ B \notin A + C $. Hence $ B + C \nparallel A + C $. 

 Since $ B + C \; \cap \; A + C = C $, and $A \ne C $, it follows that $ A \notin B + C $. Hence $ A + B \nparallel B + C $. $\Box$

\begin{prop}   
\label{ax21}
Let $ l = P + Q $ and let $ m $ be a line nonparallel to $ l $. Then either $ P \notin m $ or $ Q \notin m $. 
\end{prop}

\noindent Proof. Set $ R $ $ \equiv $ $ m \cap l $. Either $ R $ $ \ne $ $ P $ or $ R $ $ \ne $  $ Q $; let us say that $ R $ $ \ne $ $ P $. It then follows from  axiom L1 that either $ P \notin m $ or $ P \notin l $; thus $ P \notin m $. $\Box$\\

 \noindent \textit{Note.}\label{ax22} Classically, the above proposition still holds if we assume only that the line  $ m $ is  distinct from $ l $, rather than nonparallel to $ l $. This is also true constructively, as will be shown in proposition \ref{ax51}. 

\begin{prop} 
\label{ax23}
Let $ l $ be a line, let $ P $ be a point outside $ l $, 
and let $ Q $ be any point on $ l $. Then $ P + Q $ $ \nparallel $ $ l $.
\end{prop}

\noindent Proof. This follows directly from the 
definition. $ \Box $ 

\begin{prop} 
\label{lines2}
Let $l$ and $m$ be any lines.\\
\hspace*{5mm}(a) If $ \neg (l \ne m ) $, 
then $ l = m $.\\
\hspace*{5mm}(b) If $ l \ne m $, and $ n $ is  any line, then either $ n \ne l $ or  $ n \ne m $. 
\end{prop}

\noindent Proof. (a) Let $ \neg (l \ne m ) $, let 
$ P \in l $, and suppose that $ P \notin m $. Then $ l \not= m $, a contradiction; hence $ \neg (P \notin m) $, and it follows from theorem  \ref{outside2} that $ P \in m $. Thus $ P \in l $ implies $ P \in m $. Similarly, the opposite inclusion also holds.

(b) We may assume that there exists a point $ P \in l $ such that $ P \notin m $. Choose any point $ Q \in m $; thus  $ P \ne Q $.  It follows from proposition \ref{ax23} that  $ P+Q \nparallel m $. By axiom L2, either  $ n \nparallel m $ or $ n \nparallel P+Q $. In the second case, it follows from proposition \ref{ax21} that either $ P \notin n $ or $ Q \notin n $. Thus either $ n \ne l $ or  $ n \ne m $. $\Box$\\

Propositions \ref{lines1} and \ref{lines2} together show that the inequality relation ``distinct lines'' is a tight apartness.

\begin{prop}
\label{ax24}
Let $ l $  be any line, let $ R $  be a point outside $ l $, and let P and Q be distinct points on $ l $. Then $ R + P \nparallel  R + Q $. 
\end{prop}

\noindent Proof. It follows from proposition \ref{ax23} that $ R + P \nparallel l. $ It then follows from  axiom L1 that either $ Q \notin l $ or $ Q \notin R + P ;$ thus $ Q \notin R + P. $ Hence $ R + P \nparallel  R + Q.$ $\Box$

\begin{prop}
\label{ax29}
The relation ``parallel'' is an equivalence relation.
\end{prop}

\noindent Proof. Since the relation ``nonparallel'' is clearly anti-reflexive and symmetric, the relation ``parallel'' is reflexive and symmetric. Now let $ l \parallel m $ and $ m \parallel n $. Suppose that $ l \nparallel n $; then these lines have a common point. 
It follows from the parallel postulate that $ l = n $, a contradiction. This shows that $ l \parallel n $. Thus the relation ``parallel'' is transitive. $\Box$

\begin{defn}
\label{ax30}
\textnormal{An equivalence class of the relation 
``parallel''  will be called a \textit{pencil of lines.} Each pencil of lines $ \pi $ is of the form $ \pi = \pi_l \equiv \{ m \in \mathscr{L} : m \parallel l \}, $
for any line $ l $  in $ \pi $. For pencils $ \pi_l $ and 
$ \pi_m $, the expression $ \pi_l \ne \pi_m $ will mean that $ l \nparallel m $, and hence $ l' \nparallel m' $ whenever $ l' \in \pi_l $ and $ m' \in \pi_m $; we say that  \textit{``The pencils $ \pi_l $ and $ \pi_m $ are distinct.''} For any line $ l $ and any pencil $ \pi $, the expression $ l \notin \pi $ will mean that $ l \nparallel m $ for some (hence any) line $ m $ in $ \pi $.}
\end{defn}

\begin{prop}
\label{ax31}
Given any two pencils of lines, there exists a pencil distinct from each of the two given pencils. 
\end{prop}

\noindent Proof. Theorem \ref{ax20} provides three distinct pencils of lines; three applications of axiom L2 then yield the required pencil. $\Box$

\begin{prop}
\label{ax32}
Let $ l $ and $ m $ be any lines. If ~$ l = m $,
or if ~$ l \cap m = \emptyset $, then $ l \parallel m $.
\end{prop} 

\noindent Proof. This follows directly from the definition. $ \Box $ \\

The \label{ax33} statement of proposition 
\ref{ax32} is the usual classical definition of parallel lines. However, the converse is constructively invalid on the real plane $\mathbb{R}^{2}.$\footnote{See example \ref{brouB}.} 

\begin{prop}
\label{ax34}
Let $l$ and $m$ be any lines. Then $ l \parallel m $ if and only if 

\begin{displaymath}
l \cap m \ne \emptyset  ~~implies~~  l = m 
\end{displaymath}

\end{prop}

\noindent Proof. First let $ l \parallel m $, and 
let $ l \cap m \ne \emptyset $. Suppose that 
$ l \not= m $; then, by definition, $ l \nparallel m $, a contradiction. It follows from proposition \ref{lines2}(a) that $ l = m $. Now let the implication hold. Suppose that $ l \nparallel m $; thus $ l \cap m \ne \emptyset $ and $ l \not= m $, a contradiction. This shows that $ l \parallel m $. $\Box$\\

The implication in proposition \ref{ax34} is 
classically equivalent to the usual definition of parallel lines. However, the implication would not serve as a definition here; its negation is insufficient to construct a point of intersection of nonparallel lines.\footnote{See example \ref{brouI}.}

\begin{thm}
\label{ax35}
Let $l$ and $m$ be distinct parallel lines. Then $ P \ne Q $ for any point $ P $ on $  l $, and any point  $ Q $ on $ m $. Thus $ P \notin m $ for any point $ P \in l $, and $ Q \notin l $ for any point $ Q \in  m $.
\end{thm}

\noindent Proof. Let $ P $ be a point on $ l $, and let $ Q $ be a point on $ m $. We may assume that there exists a point $ R $ on $ l $ that is outside $ m $; thus $ R \ne Q $. It follows from proposition \ref{ax23} that 
$ m  \nparallel  R  +  Q $; thus $ l  \nparallel  R  +  Q $, with $ R  = l \; \cap \;  R  + Q $. Since $ Q \ne R $, it follows from axiom L1 that $ Q \notin l \,$. Thus 
$ Q \ne P $. $\Box$\\

A relation ``nonparallel''  that was not invariant with respect to the relation ``parallel'' would be unacceptable. Thus the following theorem provides one rationale for axiom L2. 

\begin{thm}
\label{ax36}
Assume for the moment only the axioms through \textnormal{L1.} Assume also that at least two distinct points lie on any given line. Then the following are equivalent:\\
\hspace*{5mm}(a) The relation ``nonparallel'' is invariant with respect to the relation ``parallel.''\\
\hspace*{5mm}(b) Axiom \textnormal{L2.}
\end{thm}

\noindent Proof. Let (a) hold, let $ l $ and $ m $ be nonparallel lines, let $ n $ be any line, and set $ P \equiv l \cap m $. Let $ n' $ be the line through  $ P $ that is parallel to $ n $. Choose a point $ Q $ on $ n' $ distinct from $ P $. It follows from axiom L1 that either $ Q \notin l $ or $ Q \notin m $; let us say $ Q \notin l $. Thus $ n' \nparallel l $ and by hypothesis it follows that $ n \nparallel l $. This proves axiom L2. The converse was proposition \ref{ax28}. $\Box$\\

 \noindent \textit{Problem.} Determine whether or not axioms L1 and L2 are independent.

\begin{lm}
\label{ax37}
Given any line $ l $, there exists a point that lies outside $ l $.
\end{lm}

\noindent Proof. Axiom G3 provides distinct points $ A $, $ B $, $ C $ with $ C \notin A + B $; it is clear that $ A + B $ $ \nparallel $ $ A + C $. By axiom L2, either $ l $ $ \nparallel $ $ A + B $ or $ l $ $ \nparallel $ $ A + C $; in the first case, set $ D $ $ \equiv $ $ l \cap  A + B $. Either $ D \ne A $ or $ D \ne B $; in the first subcase it follows from axiom L1 that either $ A \notin l $ or $ A \notin A + B $, and thus $ A \notin l $. The other three subcases are similar. $\Box$
 
\begin{thm}
\label{ax38}
At least two distinct points lie on any given line.
\end{thm}  

\noindent Proof. Let $ l $ be any line. Use lemma \ref{ax37} to construct a point $ R $  that is 
outside $ l $. Let three non-collinear points be given as in axiom G3, and construct the three lines formed by these points. By theorem \ref{ax20} these lines are nonparallel in pairs. It follows from two applications of axiom L2 that two of the three lines are nonparallel to $ l $; denote them $ m $ and $ n $.

 Let $ m' $ and $ n' $ be the lines through  $ R $ such that $ m' \parallel m $ and $ n' \parallel n\,$. Thus $ m' \nparallel n' $, $ l \nparallel m' $, $ l \nparallel n' $, and $ R  =  m' \cap n' $. Set $ P \equiv l \cap m' $, and set $ Q \equiv l \cap n' $. Since $ R \notin l $, we have $ R \ne Q $. Thus it follows from axiom L1 that either $ Q \notin n' $ or $ Q \notin m' $. Thus   
$ Q \notin m', $ and it follows that $ Q \ne P $. $\Box$

\begin{cor}
\label{ax39}
Any line may be expressed in the form $ l = P + Q $. 
\end{cor}

\begin{cor}
\label{ax40}
Given a line $ l $, and any point $ P $ on $ l $, there exists a point $ Q $ on $ l $ that is distinct from $ P $. 
\end{cor}

Theorem \ref{ax38} now enables the proof of the next proposition, which was foretold in the note following proposition \ref{ax21}. 

\begin{prop}
\label{ax51} 
Let $ l = P + Q $ and let $ m $ be a line distinct from $ l $. Then either $ P \notin m $ or $ Q \notin m $. 
\end{prop}

\noindent Proof. Either there exists a point $ R \in m $ that is outside $ l $, or there exists a point $ S \in l $ that is outside $ m $.

In the first case, $ R + P \nparallel R + Q $, and it follows from axiom L2 that $ m $ is nonparallel to one of these lines; we may assume that $ m \nparallel R + P $. Since $ P \ne R = m \cap R + P $, it follows from axiom L1 that  $ P \notin m $.

 In the second case, use corollary  \ref{ax39} to express $ m $ in the form $ m = U + V $. Using $ S $ in place of $ R $, it follows from the first case, with the lines reversed,  that either $ U \notin l $ or $ V \notin l $; let us say that $ U \notin l $. Now, using the point $ U $ in place of $ R $, the first case  shows that either $ P \notin m $ or $ Q \notin m $. $\Box$\\

When distinct lines are given, the following corollary circumvents the need to consider both alternatives in the definition. 

\begin{cor}  
\label{ax52}
Let $ l_1 $ and $ l_2 $ be distinct lines. If $ l $ is either of these lines, then there exists a point on $ l $ that is outside the other line. 
\end{cor}

\begin{thm}
\label{ax41}
There exist bijections mapping

 (a) the lines in any two pencils of lines,

 (b) the points on any two lines,

 (c) the points on any line and the lines in any pencil of lines.
\end{thm} 

\noindent Proof. First let $ \pi $ be any pencil of lines, and let $ l $ be any line with  $ l \notin \pi $. Define a map $ \psi : \pi \rightarrow l $ as follows. For any line $ m $ in $  \pi $,  set $ \psi (m)  \equiv  m \cap l $. 
Let $ m $ and $ n $ be lines in $ \pi $ with $ m  \ne  n $. Set $ P \equiv \psi (m) = m \cap l $, and $ Q \equiv \psi (n) = n \cap l $. It follows from  theorem \ref{ax35} that  $ P \ne Q $. This shows that the map $ \psi $ is injective; it is onto $ l $ because of the parallel postulate. If $ P  = m \cap l $ and $ Q  = n \cap l $ are distinct points on $ l $, then it follows from axiom L1 that $ Q \notin m $; thus $ m \ne n $. This shows that the inverse  is injective, and thus the map $ \psi : \pi \rightarrow l $ is a bijection. 

Now the required maps are obtained by combining various instances of the map $ \psi $ and its inverse.
$\Box$\\

 \noindent \textbf{Appendix A to section \ref{sec:ax}; Alternative axioms.} The axioms in group \textbf{L} are classically equivalent to the uniqueness portions of axioms G1 and G2. Thus axiom group \textbf{L} could be deleted, with group \textbf{G}  rewritten as shown below. This would have the advantage of showing clearly that no essentially new axioms are introduced; the traditional axioms are only rewritten as strong classical equivalents. On the other hand, the procedure followed above in section \ref{sec:ax} has the following advantages: (i) The axioms in group \textbf{G} retain their traditional form. (ii) Axiom group \textbf{L} clearly indicates what supplements are needed for the constructivization. \\

\noindent \textbf{Axiom G1*.} \textit{If $ P $ and $ Q $ are distinct points, then there exists a line $ l $ that passes through both points $ P $ and $ Q $. The line $ l $ is unique in the following sense: If $ l_1 $ and $ l_2 $ are distinct lines, both through $ P $, then either $ Q $ lies outside $ l_1 $, or $ Q $ lies outside $ l_2 $.}\\

\noindent \textbf{Axiom G2*.} \textit{If $ P $ is any point and $ l $ is any line, then there exists a line $ m $ that passes through $ P $ and is parallel to $ l $. The line $ m $ is unique in the following sense: If $ m_1 $ and $ m_2 $ are distinct lines, both through $ P $, then either $ m_1 $ is nonparallel to $ l $, or $ m_2 $ is nonparallel to $ l $.}\\

\noindent \textbf{Appendix B to section \ref{sec:ax}; 
Heyting's axiom system.}  Arend Heyting has introduced axioms for incidence geometry [He2], with the goal of extending the resulting plane to a projective plane. The express purpose of Heyting's axiom system, and the divergent approaches of [He2] versus the present paper, make direct comparisons difficult, if not meaningless. Nevertheless, a few comments on the various features of the two systems might not be out of place. A cursory comparison, rather than an exhaustive discussion, is intended. 

Heyting uses unfamiliar notation for certain concepts. His reason, it may be surmised, was to provide clear indicia to remind the reader that intuitionistic mathematics is different than classical mathematics; the unusual notation was perhaps meant to serve as a reminder that classical ideas, such as the \textit{principle of the excluded middle,} were not to be invoked. In particular, the usual symbol $\ne$ for denoting \textit{distinct points} and \textit{distinct lines}, and the symbol $\notin$ for \textit{point outside a line}, were eschewed; perhaps this was meant to emphasize the special meanings adopted for these concepts, and avoid the possibility of conjuring up the idea of simple negation. 

In contrast, at a later time, Bishop strived to demonstrate that constructive mathematics was not a new type of mathematics, but rather a return to older, more stringent standards. Thus Bishop tried to make constructive mathematics look the same as traditional mathematics; this is the approach followed in the present paper. 

The axiom system in [He2] involves quite a few more axioms than the three axioms of the traditional system, and many of the axioms lack an immediately clear intuitive interpretation. In contrast, guided by the idea of the preceding paragraph, the present paper attempts to obtain an axiom system similar to the classical system; it is identical to the traditional system in the sense outlined in Appendix A above. 

Both systems make use of apartness relations, which are natural in constructive mathematics.\footnote{The terminology for apartness relations is somewhat variegated in [He2] and [MRR]; the terminology of the latter is used in the present paper.} The apartness concept was developed by A. Heyting (1898-1980) in his 1925 dissertation, written under the direction of L. E. J. Brouwer. 

In [He2], the parallel postulate is not taken as a simple axiom, but proved later, as [He2, Theorem 3], using the more complicated axioms. The present paper uses the traditional axiom G2. 

Axiom A4 in [He2] follows directly from proposition \ref{ax21} in the present paper. 

Axiom A5 in [He2] follows from theorem \ref{ax35} in the present paper. 

Axiom A7(ii) in [He2] requires that every line contain at least four points; this excludes the four-point and nine-point geometries.\footnote{The possibility of the four-point geometry sprang up as one of the hurdles in the present work; see theorem \ref{des9}, condition (a) therein, and the subsequent note labeled \textit{Problems}.} In axiom G3, the present paper assumes only the existence of three non-collinear points, as is traditional; it is then proved in theorem \ref{ax38} that at least two distinct points lie on any given line. 

In Section 4 (Definition 2) of [He2], the definition of \textit{distinct lines} is not symmetric as stated; the concept is later shown to be symmetric (Theorem 2). In the present paper, definition \ref{ax4} for \textit{distinct lines} is symmetric, and weaker than the definition in [He2]; the concept is later shown, in corollary \ref{ax52}, to have the stronger property. 

In [He2], the notion of \textit{parallel lines,} given in Definition 5, does not include \textit{equal} lines; \textit{parallelism} is not an equivalence relation. The definition in the present paper more closely matches the traditional idea. 

The terms \textit{parallel lines} and \textit{nonparallel lines} (the later being termed ``intersecting lines'' in [He2]) are given separate definitions in [He2]. In the present paper, \textit{nonparallel} is given the primary, affirmative definition; \textit{parallel} is the negation of \textit{nonparallel}, and propositions \ref{ax32} and \ref{ax34} connect the concept of \textit{parallel lines} with the traditional definition, as far as is constructively possible.  

The above comments notwithstanding, one must remember that the axiom system in [He2] was written expressly to enable the extension to a projective plane. It remains an open problem to use the axiom system and methods of the present paper to effect this extension. (The sole concern of the present paper is the coordinatization of the plane.) 

\textit{Acknowledgement.} Many thanks are due the referee for bringing Heyting's paper to the author's attention, and for suggesting the addition of this appendix.

\section{Dilatations}
\label{sec:dil}

The \label{dila1} symmetries required for Desarguesian geometry are mappings of the plane that preserve direction. On the real plane $\mathbb{R}^{2},$ these include uniform motions and expansions about a fixed point. 

\begin{defn}
\label{dila2}
\textnormal{A \textit{dilatation} is a map 
$ \sigma : \mathscr{P} \rightarrow \mathscr{P} $  
such that:\\
\hspace*{5mm}(a) $ \sigma $ is injective and onto, and \\
\hspace*{5mm}(b) if $ P $ and $ Q $ are distinct points, then 
$\; P + Q \; \parallel \; \sigma P + \sigma Q.  $}
\end{defn} 

The normally-used classical  condition, weakly injective, would not suffice here. The  \textit{direction-preserving property}\label{dirprop} \ref{dila2}(b) cannot even be stated unless the map is injective.  

Classically \label{dila5}[A], a ``dilatation'' need not be even weakly injective; if not, it is termed ``degenerate.'' Although only nondegenerate dilatations are of use, crucial steps in the most important classical constructions depend on the notion that a dilatation must be either degenerate or injective; this notion is constructively invalid on the real plane $\mathbb{R}^{2}$.\footnote{See example \ref{brouL}.}\\

The following form of the direction-preserving  property is often convenient.

\begin{lm} 
\label{dila6}
Let $ \sigma $ be a dilatation, let $ P $ and $ Q $ be distinct points, and let $ l $ be a line parallel 
to $ P + Q $. If $ \sigma P \in l ,$ then 
also $ \sigma Q \in l $.
\end{lm}

In the construction of the division ring in section \ref{sec:ring}, we will need to construct a dilatation in a situation where the onto property is not immediate; the following  theorem will be required. 

\begin{thm}
\label{dila7}
Let $ \sigma $ : $\mathscr{P} \rightarrow \mathscr{P}$ be a map such that \\
\hspace*{5mm}(a) $ \sigma $ is injective, and \\
\hspace*{5mm}(b) if $ P $ and $ Q $ are distinct points, then
$ P + Q \; \parallel \; \sigma P + \sigma Q $.\\ 
 \noindent Then $ \sigma $ is a dilatation.
\end{thm} 

\noindent Proof. We need only prove that $ \sigma $ is onto. Let $ S $ be any point. Using axiom G3 and 
theorem \ref{ax20}, construct  three distinct noncollinear 
points $ A ,$ $ B ,$ $ C $,  forming three 
nonparallel lines. It follows that the points 
$ \sigma A ,$ $ \sigma B ,$ $ \sigma C $ have the same properties. Applying axiom L1, we may choose two of these points, which we may denote 
by $ \sigma P $ and $ \sigma Q ,$ such that $ S \notin \sigma P + \sigma Q $.

Set  $ l \equiv \sigma P + \sigma Q $ and $ l' \equiv P + Q $; thus $ l \parallel l' $. Set $ m \equiv \sigma P + S $ and $ n \equiv \sigma Q + S $. It follows from  proposition \ref{ax24} that $ m \nparallel n $; thus $ S = m \cap n $. Let $ m' $ and $ n' $ be the lines such 
that $ P \in m' \parallel m ,$ and $ Q \in n' \parallel n $.
Thus $ m' \nparallel n' $; set $ R \equiv m' \cap n' $.

Either $ R \ne P $ or $ R \ne Q $. In the first case, $ R \ne m' \cap l' $; it follows from axiom L1 that $ R \notin l' $, and thus $ R \ne Q $. The second case is similar; thus in either case we have both $ R \ne P $ and $R \ne Q $. 

Since $ \sigma P \in m \parallel m' = P + R, $ it follows from lemma \ref{dila6} that 
$ \sigma R \in m $. Similarly, $ \sigma R \in n $; hence $ \sigma R = S $. Thus $ \sigma $ is onto. $\Box$\\

Although the inverse of a dilatation is weakly injective,  the stronger condition, injective, is not evident. It will be required that dilatations have inverses that are also dilatations; thus the inverses must be injective. One might consider specifying this in the definition, but this would create serious difficulties in constructing dilatations. The next theorem settles this problem. 

\begin{thm}
\label{dila9} 
If $\sigma$ is a dilatation, then the inverse $ \sigma ^{-1} $ is also a dilation.
\end{thm}

\noindent Proof. We will first show that $ \sigma^{-1} $ is injective. Since $\sigma$ is injective and onto, any point may be expressed uniquely
in the form $ R'  \equiv \sigma R $.

Let $ S' $ and $ T' $ be points with $ S' \ne T' $, and
set $ w'  \equiv  S' + T' $. Choose any line $ l' $ with $ l' \nparallel  w' $, and set $ V' \equiv l' \cap w' $.
Let $ l $ and $ w $ be the lines through  $ V $ that are 
parallel to $ l' $ and $ w' $. Using theorem \ref{ax38}, construct a point $ P $ on $ l $ distinct from $ V $. Thus $ P' \in l' ,$ and $ P' \ne V' $; it follows from axiom L1 that $ P \notin  w $ and $ P' \notin  w' $.

Set $ m' \equiv P' + S' $ and
$n' \equiv P' + T' $;
thus $m' \nparallel   n' ,$
~$ m'  \nparallel  w',$ and ~$ n' \nparallel w'$.
Let $ m $ and $ n $ be the lines through $ P $ that are parallel to
$ m' $ and $ n' $; set
$A \equiv  m  \cap w $ and set $B \equiv  n  \cap w $. Since $ P \notin w ,$ we have
$A \ne P  =  m  \cap  n $;
it follows that $A \notin   n $.
Thus $A \ne B$.

Either $V \ne A$ or $V \ne B$.
In the first case, 
$V' \in w' \parallel w = V + A ,$
and it follows from lemma \ref{dila6} that $A' \in  w' $. 
Now $A' \in w' \parallel w = A + B ,$
and it follows that $B' \in w'$.
This shows that  both $A'$ and $B'$ lie on $w'$;
the second case produces the same result.

Since 
$ P' \in m' \parallel m = P + A ,$
we have $A' \in  m' $. Similarly, $B' \in n' $.
Thus $A' = m' \cap w' = S' $, and $ A = S $. Similarly, $B' = T'$, and $ B = T $. 
Since $A \ne B$, this means that  $ S \ne T$.
This shows that $\sigma^{-1}$ is injective.

Now let $ P' $ and $ Q' $ be any distinct points.
Then $ P \ne Q ,$ and  
$\sigma^{-1} P' + \sigma^{-1} Q' = P + Q  \parallel 
\sigma P + \sigma Q = P' + Q'$. Hence $\sigma^{-1} $ is a dilatation. $\Box$

\begin{thm}
\label{dila10} 
The dilatations form a group D.
\end{thm}

\noindent Proof. The identity $ 1 : \mathscr{P}  \rightarrow \mathscr{P}  $ is clearly a dilatation. Theorem \ref{dila9} has shown that the inverse of a dilatation is also a dilatation. Let $ \sigma_1 $ and $ \sigma_2 $ be dilatations; the product is clearly injective and onto. If $ P $ and $ Q $ are distinct points, then
$ P + Q \parallel \sigma_2 P + \sigma_2 Q \parallel \sigma_1 \sigma_2 P + \sigma_1 \sigma_2 Q $. Hence $ \sigma_1 \sigma_2 $ is a 
dilatation. $\Box$

\begin{defn}
\label{dila11}
\textnormal{Let $ \sigma $ be a dilatation.
A point $ P $ is a \textit{fixed point of} $ \sigma $ if $ \sigma P = P $. If $ \neg ( \sigma P = P ) $ for all points $ P ,$ then $ \sigma $ \textit{has no fixed point.}}\footnote{The definition of ``has no fixed point'' is uncharacteristically negativistic for a constructive theory. It is given mainly to enable a discussion of the definition of ``translation'' in section \ref{sec:tran}. In contrast, theorem \ref{tran9} will provide a strong version of this idea.} 
\end{defn}

\begin{lm}
\label{dila12}
Let $ \sigma $ be a dilatation with two distinct fixed points $ P $ and $ Q $. If $ R $ is any point outside the line $ P + Q ,$ then $ \sigma R = R $.
\end{lm}

\noindent Proof. Set $ l \equiv P + Q , $ set $  m \equiv P + R ,$ and set $ n \equiv Q + R $.
Since $ \sigma P \in m = P + R ,$ it follows from
lemma \ref{dila6}  
that $ \sigma R \in m $.
Similarly, $ \sigma R \in n $.
Thus $\sigma R = m \cap n = R$. $\Box$

\begin{lm}
\label{dila13}
If $ \sigma $ is a dilatation with two distinct fixed points, then $ \sigma = 1 $. 
\end{lm}

\noindent Proof. Let $ P $ and $ Q $ be distinct fixed points, and 
set $ l \equiv P + Q $. Using Lemma 
\ref{ax37}, construct any point S outside $ l $. It follows from  
lemma \ref{dila12} that 
$ S $ is a fixed point. 
Since $ S \ne P $ and $ S \ne Q ,$ 
the lemma applies also to the lines 
$ m \equiv P + S $ and $ n \equiv Q + S $. 
It also follows that   
$ l \nparallel m ,$ 
$ l \nparallel n ,$  
$ P = l \cap m ,$ and $ Q = l \cap n $.
 
Now consider any point $ R $. 
Either $ R \ne P $ or $ R \ne Q $; we may 
assume the first case. Thus $ R \notin l $ 
or $ R \notin m $; in either case, it follows from lemma \ref{dila12} that $ \sigma R = R $ . 
This shows that $ \sigma = 1 $. $\Box$

\begin{thm} \label{dila14}  
A dilatation is uniquely determined by the images of any two distinct points.
\end{thm}

\begin{prop}
\label{dila15}
Let $ \sigma_1 $ and $ \sigma_2 $ be dilatations. If $ \sigma_1 \sigma_2 \ne 1 ,$ then either $ \sigma_1 \ne 1 $ or $ \sigma_2 \ne 1 $. 
\end{prop}

\noindent Proof. We have $ P \ne \sigma_1 \sigma_2 P $ 
for some point $ P $. Since $ \sigma_1^{-1} $  is injective, we 
have  $ \sigma_1^{-1} P \ne \sigma_2 P $.
Either $ P \ne \sigma_1^{-1} P$ or 
$ P \ne \sigma_2 P $.
In the first case, $ \sigma_1 P \ne P ,$
and thus $ \sigma_1 \ne 1 $.
In the second case, $ \sigma_2 \ne 1$. $\Box$

\begin{defn}
\label{dila16}  
\textnormal{A \textit{trace} of a 
dilatation $ \sigma $ is any line of the 
form $ P + \sigma P ,$ 
where $ P \ne \sigma P $. 
For any dilatation $ \sigma ,$ 
the set of lines} 
\begin{displaymath}
t( \sigma ) \equiv \{ l \in \mathscr{L} : 
l \; \textnormal{is a trace of} \; \sigma \} 
\end{displaymath}
\textnormal{will be called the 
\textit{trace family of} $ \sigma $.}
\end{defn}

\begin{prop}
\label{dila17}
Let $ \sigma $ be a dilatation. 
Then \\
\hspace*{5mm}(a) $ t( \sigma ) = \emptyset $ if and only if $ \sigma = 1 ,$\\
\hspace*{5mm}(b) $ t( \sigma) \ne \emptyset $ if and only if $ \sigma \ne 1 ,$\\
\hspace*{5mm}(c) 
$ t( \sigma^{-1}) = t( \sigma )$.
\end{prop}

\label{dila18}
In the classical theory [A], the definition of ``trace'' is slightly weaker, with the result that $t(1) = \mathscr{L} $; this is of little consequence, since it is presumed that one always knows whether $ \sigma = 1$ or $ \sigma \ne 1$. For a constructive development, our definition is more convenient: if a dilatation $\sigma $ has a trace, then $ \sigma \ne  1 $. In general, we will not know whether $ \sigma = 1 $ or $ \sigma \ne 1 ,$ or whether $ t(\sigma) = \emptyset $ or $ t (\sigma) \ne \emptyset $.\footnote{See example \ref{brouE}.}

\begin{lm}
\label{dila19}
If $ l $ is a trace of a 
dilatation $ \sigma ,$ then there exist at least two distinct points $ P $ and $ R $ such that
   $ l = P + \sigma P $
and $ l = R + \sigma R $. 
\end{lm}

\noindent Proof. Choose  a point $ P $ such that 
$ l = P + \sigma P $ and set $ R \equiv \sigma P$.
Since $ \sigma P \in l = P + R ,$ it follows from
lemma \ref{dila6}  that $ \sigma R \in l$.
Since $ R \ne P,$ we have 
$ \sigma R \ne \sigma P = R $;
thus $ \l = R + \sigma R $. $\Box$

\begin{thm}
\label{dila20}
Let $ \sigma $ be a dilatation, and let $ l $ be a trace of $ \sigma $. If $ Q $ is a point on $ l ,$ then $ \sigma Q $ also lies on $ l $.
\end{thm}

\noindent Proof. Using lemma \ref{dila19}, choose a 
point $ P $ such that
$ l = P + \sigma P $
and $ Q \ne P $.
Since $ \sigma P \in \l = P + Q ,$
it follows from lemma 
\ref{dila6} that
$ \sigma Q \in \l $. $\Box$

\begin{cor}
\label{dila21}
Let $ \sigma $ be a dilatation. 
The intersection of any two nonparallel traces 
of ~$ \sigma $ is a fixed point.
\end{cor}

\begin{lm}
\label{dila22}
Let $ \sigma $ be a dilatation with fixed point $ P $. Then every trace of $ \sigma $ passes through $ P $.
\end{lm}

\noindent Proof. Let $ l $ be a trace of $ \sigma $. Using lemma \ref{dila19}, choose a point $ Q $ such that $ l = Q + \sigma Q$ and $ P \ne Q $.
Since $ P + Q \parallel \sigma P + \sigma Q = P + \sigma Q ,$ it follows that $ \sigma Q \in P + Q $.
Since the line $ P + Q $ passes through both points $ Q $ and $ \sigma Q$, we have $ P + Q = l $. $\Box$

\begin{lm}
\label{dila23}
Let $ \sigma $ be a dilatation with fixed point $ P ,$
and let $ l$ be a line through $ P $. If $ Q $ is a point on $ l ,$ then $ \sigma Q $ also lies on $ l $.
\end{lm}

\noindent Proof. Suppose that $ \sigma Q \notin l $. Then $ \sigma  Q \ne P$; applying $ \sigma^{-1},$ we have $ Q \ne P $. Thus $ l = P + Q $. Also, $ \sigma Q \ne Q $;  thus $ Q + \sigma Q $ is a trace of $ \sigma $. It follows from  lemma \ref{dila22} that  
$ P \in Q + \sigma Q $; thus $ Q + \sigma Q = l ,$ and
$ \sigma Q \in l ,$ a contradiction. Hence $ \sigma Q \in l $. $\Box$\\

If $ \sigma $ is a dilatation with $ \sigma \ne 1 ,$ then it follows from lemma \ref{dila13} that $ \sigma $ may have at most one fixed point. More precisely, 
if $ P $ and $ Q $ are fixed points, then $ P = Q $. To establish theorem \ref{dila25}, we will require the following stronger result: 

\begin{lm}
\label{dila24} 
Let $ \sigma \ne 1 $ be a dilatation with fixed point $ P $. If $ Q $ is any point distinct from $  P, $ then $ \sigma Q \ne Q $.
\end{lm}

\noindent Proof. Since  $ \sigma \ne 1 $, we may construct at least one point $ R $ such that $ \sigma R \ne R .$ Set $ l \equiv R + \sigma R $; it follows from  lemma \ref{dila22} that  $ P \in l $. 

Let $ S $ be any point outside $ l $. Set $ m \equiv R + S $, and $ m' \equiv \sigma R + \sigma S $; thus $ m \parallel m' $. Since $ S \notin l $, we have $ m \nparallel l $. Since $ \sigma R \ne R = l \cap m $, it follows that $ \sigma R \notin m $. Thus 
$ m \ne m' $, and it follows from theorem \ref{ax35} that $ \sigma S \ne S $. 

Construct one point $ T $ outside $ l $, and set $ l' \equiv T + \sigma T $. It follows that $ P \in l' $, and $  l  \nparallel  l' $. Also, for any point $ S $ outside $ l' $ we have $ \sigma S \ne S $. Now let $ Q $ be any point distinct from $ P $. It follows from  axiom L1 that either $ Q \notin l $ or $ Q \notin l' $; hence $ \sigma Q \ne Q $. $\Box$

\begin{thm}
\label{dila25}
Let $ \sigma \ne 1 $ be a dilatation with fixed point $ P $. Then

\begin{displaymath}
t(\sigma) = \mathscr{L}_P \equiv \{ l \in \mathscr{L} : P \in l \} 
\end{displaymath} 
\end{thm}

\noindent Proof. Lemma \ref{dila22} shows that $ t ( \sigma) \subseteq \mathscr{L}_P $. Now let $ l $ be any line through $ P $. Construct a point $ Q $ on $ l $ such that $ Q \ne P  $; it follows from lemma \ref{dila23} that $ \sigma Q \in l $, and from lemma \ref{dila24} that $ \sigma Q \ne  Q $. Thus $ l = Q + \sigma Q ,$ and this is a trace of $ \sigma $. $\Box$

\section{Translations}
\label{sec:tran}

The translations of the geometry are the symmetries which we perceive as uniform motions. These maps form the substructure of the coordinatization to be carried out in sections \ref{sec:ring} and \ref{sec:coor}. Points on the plane will be located using translations. Selected translations will determine the unit points on the axes. The scalars will be  certain homomorphisms of the translation group; they will relate the various translations to each other, and provide coordinates for the points.

Given a dilatation $ \tau ,$ the classical definition \label{tran1} of a translation, ``Either $ \tau = 1 ,$ or $ \tau $ has no fixed point,'' is constructively invalid for translations on the real plane $ \mathbb{R}^{2} $.\footnote{See example \ref{brouE}.}
From the following list of classically equivalent conditions, we give a Brouwerian counterexample for the first, prove the equivalence of the last three, and choose the last for a definition.\footnote{Condition (d) is classically equivalent to conditions (b) and (c) only under the concept of trace adopted in definition \ref{dila16}.}\\
\hspace*{5mm}(a) Either $ \tau = 1 ,$ or $ \tau $ has no fixed point.\\
\hspace*{5mm}(b) If $ \tau $ has a fixed point, then $ \tau = 1 $.\\
\hspace*{5mm}(c) If $ \tau \ne 1 ,$ then $ \tau P \ne P $ for all points $ P $.\\
\hspace*{5mm}(d) Any traces of $ \tau $ are parallel.\\
Condition (d) is chosen for its simplicity and intuitive imagery. All of the last three conditions will be required; their equivalence will be demonstrated in theorems \ref{tran3} and \ref{tran9}. 

\begin{defn}\label{tran2} 
\textnormal{A dilatation $\tau$ will be called a \textit{translation} if any traces of $\tau$ are parallel.}
\end{defn}

\begin{thm}\label{tran3}
A dilatation $\tau$ is a translation if and only if the following implication is valid: 
If $ \tau $ has a fixed point, then $ \tau = 1 $.
\end{thm}

\noindent Proof. First, let $\tau$ be a translation and let $ P $ be a fixed point of $\tau$. Suppose that $\tau \ne 1$. It then follows from theorem \ref{dila25} 
that the trace family $ t(\tau) $ is the family $ \mathscr{L}_P ,$ which contains nonparallel lines; this is a contradiction. Hence $ \tau = 1 $.

Now assume the implication and let $ l $ and $ m $ be any traces of $ \tau $. Suppose that $ l \nparallel m $; it then follows from corollary \ref{dila21} that the point $ P \equiv l \cap m $ is a fixed point. Thus $\tau = 1 ,$  and $\tau$ has no traces, a contradiction. Hence $ l \parallel m $. $\Box$

\begin{cor}\label{tran4}
Let $ \tau $ be a dilatation. If either $ \tau = 1 $ or $\tau$ has no fixed point, then $\tau$ is a translation.
\end{cor}

Although the condition of   corollary \ref{tran4} is the classical definition of ``translation,'' the converse is constructively invalid on the real 
plane $ \mathbb{R}^{2} $.\footnote{See example  \ref{brouE}.}  

\begin{lm}\label{tran6}
Let $\tau \ne 1 $ be a translation, and let $ l $ be a trace of $ \tau $. Then $ \tau Q \ne Q $ for any point $ Q $ outside $ l $.
\end{lm}

\noindent Proof. Choose a point $ P $ such that 
$ l = P + \tau P $. Since $ Q \notin l ,$ 
we have $ P + Q \nparallel P + \tau P $.
It follows from  axiom L1 that $ \tau P \notin P + Q $.
Thus the lines $ P + Q $ 
and  $ \tau P + \tau Q $ are parallel and distinct. It now follows from theorem 
\ref{ax35} 
that $\tau Q \ne Q$. $\Box$

\begin{lm}
\label{tran7}
Let the lines $ l $ and $ m $ be parallel and distinct. For any point $ Q ,$  either $ Q \notin l $ or $ Q \notin m $. 
\end{lm}

\noindent Proof. Choose a line $ n $ containing $ Q $ such that $ n \nparallel l $; thus also   $ n \nparallel m $. Set $ P \equiv n\cap l $, and $ R \equiv  n\cap m $. It follows from  theorem \ref{ax35} that $  P \ne R $; we may assume that $ Q \ne P$. It now follows from axiom L1 that $ Q \notin l $. $\Box$\\

Although it follows from theorem \ref{tran3} that a translation $ \ne 1 $ has no fixed point, we will require the following  stronger result.

\begin{thm}
\label{tran9}
Let $\tau $ be a translation. If $\tau \ne 1 ,$ then $ \tau P \ne P $ for every point $ P $. 
\end{thm}

\noindent Proof. Since $ \tau \ne 1 $ we may choose at least one point $ Q $ such that $ \tau Q \ne Q $.
Set $ l \equiv Q + \tau Q , $ and choose any point $ R $ outside $ l $. It follows from  lemma \ref{tran6} that
 $ \tau R \ne R $; set $ m \equiv R + \tau R $. Since $ l $ and $ m $ are traces of $ \tau , $  they are parallel; they are also distinct. 

Now let $ P $ be any point; it follows from lemma \ref{tran7} that either $ P \notin l $ or $ P \notin m $. In either case, lemma \ref{tran6} shows that $ \tau P \ne P $. $\Box$\\

For translations, the next lemma extends theorem \ref{dila20}. 

\begin{lm}
\label{tran12}
Let $\tau $ be a translation,  let $\pi $ be a pencil of lines with $ t ( \tau) \subseteq \pi $, and let  $l $ be any line in $  \pi $. If $ P \in l, $  then also $ \tau P \in l $.
\end{lm}

\noindent Proof. Suppose that $ \tau P \notin l $. Then 
$\tau P \ne P, $  and $ m \equiv P + \tau P$ is a trace of $\tau $. Since both $ l $ and $ m $ are in the pencil $ \pi, $ they are parallel, and thus equal, a contradiction. This shows that $ \tau P \in l $. $\Box$\\

The statement: ``For any translation $ \tau ,$ there exists a 
pencil of lines $ \pi $ such that 
$ t ( \tau ) \subseteq \pi $'' is constructively 
invalid.\footnote{See example \ref{brouK}.}  

\begin{thm}
\label{tran14}
The traces of a translation 
$ \tau \ne 1 $ form a pencil of lines.
\end{thm}

\noindent Proof. Choose any trace $ l $ of $ \tau $.  
Since any traces of $ \tau $ are parallel, it follows that 
$ t ( \tau ) \subseteq \pi_l $.  
Now let $ m \in \pi_l , $  and choose any 
point $ Q $ on $ m $. Lemma
\ref{tran12} shows that $ \tau Q \in m ,$ and it follows from  
theorem \ref{tran9} that
$ \tau Q \ne Q$. Thus $ m = Q + \tau Q , $  
and $ m $ is a trace of $ \tau $. $\Box$

\begin{defn}
\label{tran15}
\textnormal{The trace pencil of a translation $ \tau \ne 1 $ will be called the \textit{direction of} $ \tau $.}
\end{defn}

\begin{thm}
\label{tran17}
A translation is uniquely determined by the image of a single point. 
\end{thm}

\noindent Proof. Let $ \tau_1 $ and $ \tau_2 $ be translations such that $ \tau _1 P = \tau _2 P $ for some point $ P $.

Consider first the special case in which
$ \tau _1 P = \tau _2 P \ne P $. Denote by $ \tau $ either of the given translations. Thus $ \tau \ne 1 $, and $ l \equiv P + \tau P $ is a trace of $ \tau $.
Select a point $ Q $ outside the line $ l $; thus $ P \ne Q $. Let $ l' $ be the line through $ Q $, parallel to $ l $. Then $ l' $ is also a trace of $ \tau ,$ 
and it follows from theorem \ref{dila20} that $ \tau Q \in l' $. 

Set $ m \equiv P + Q $, and let $ m' $ denote the line parallel to $ m $ through $ \tau P $. It follows from Lemma \ref{dila6} that $ \tau Q \in m' $. Since $ l \nparallel m $, also $ l' \nparallel m' $; thus $ \tau Q = l' \cap m' $. The lines $ l, ~ l', ~m, ~m' $ are uniquely determined, solely by $ P, ~\tau P, $ and $ Q $. Thus the point $ \tau Q $ is uniquely determined.

Since the point $ \tau Q $ was determined independently of the choice $ \tau = \tau _1 $ or $ \tau = \tau _2  ,$  we have $ \tau _1 Q = \tau _2 Q $. Since $ \tau _1 $ and $ \tau _2 $ agree at two distinct points, it follows from theorem \ref{dila14} that $ \tau _1 = \tau _2 $.

Now consider the general case, and suppose that $ \tau _1 \ne \tau _2 $. Suppose further that $ \tau _2 P \ne P $.  Using the special case, we have $  \tau _1 = \tau _2 ,$ a contradiction; hence $ \tau _2 P = P $. From theorem \ref{tran3} it follows that $ \tau _2 = 1 ,$  and also  that $ \tau _1 = 1 $, a contradiction. This shows that $ \tau _1 = \tau _2 $. $\Box$

\begin{thm} \hspace*{5mm}
\label{tran18} 

(a) The translations form a 
group $ T $.

(b) $ T $ is an invariant subgroup of the dilatation group $ D $.

(c) Let ~$ \tau $ be a translation, and 
let ~$ \sigma $ be a dilatation. Then
~$ t ( \sigma \tau \sigma ^ {-1} ) = 
t ( \tau ) $.
\end{thm}

\noindent Proof. (a) It follows directly from theorem 
\ref{tran3} that $ 1 \in T, $  
and that $ \tau ^{-1} \in T $ whenever 
$  \tau \in T $.  
Now let $ \tau _1, \tau _2 \in T ,$ and let
$ \tau _1 \tau _2 P = P $
for some point $ P $.
Then $ \tau _2 P = \tau _1^{-1} P $. 
It follows from theorem \ref{tran17} 
that $ \tau _2 = \tau _1^{-1},$   
and thus 
$ \tau _1 \tau _2 = 1 $. 
This shows that $ \tau _1 \tau _2 \in T $.

(b) Let $ \tau \in T $, let  
$ \sigma \in D , $  and let
$ \sigma \tau \sigma ^{-1} P = P $ for some point P. Then
$ \tau \sigma ^{-1} P = \sigma ^{-1} P $; it follows 
that $ \tau = 1 , $  
and thus $ \sigma \tau \sigma ^{-1} = 1 $.  
This shows that 
$ \sigma \tau \sigma ^{-1} \in T $.

(c) Let $ l = Q + \sigma \tau \sigma ^{-1} Q $ be a trace of $ \sigma \tau \sigma ^{-1},$ and set $ P \equiv \sigma ^{-1} Q $. Then $ l = Q + \sigma \tau \sigma ^{-1} Q \parallel \sigma ^{-1} Q + \tau \sigma ^{-1} Q = P + \tau P $, a trace of $ \tau $. It follows from theorem \ref{tran14} that $ l $ is also a trace of $ \tau $. This shows that $ t(\sigma \tau \sigma^{-1}) \subseteq t(\tau) $ for all $ \tau \in T $ and all $ \sigma \in D $. In this inclusion replace $ \sigma $ by $ \sigma ^{-1} $ and then $ \tau $ by $ \sigma \tau \sigma ^{-1} $; thus $ t ( \tau ) \subseteq t ( \sigma \tau \sigma ^{-1} ) $. $\Box$

\begin{thm}
\label{tran19}
Let $ \pi $ be any pencil of lines. Then 
\begin{displaymath} 
T_\pi \equiv \{ \tau \in T : \hspace*{1mm} t(\tau) \subseteq \pi \} 
\end{displaymath} 
is a subgroup of $ T $.
\end{thm}

\noindent Proof. It is clear from proposition 
\ref{dila17} that $ 1 \in T_\pi ,$ and that $ \tau ^{-1} \in T_\pi $ whenever $ \tau \in T_\pi $. Now let 
$ \tau _1, ~\tau _2 \in T $, and let $ l = P + \tau _1 \tau _2 P $ be a trace of $  \tau _1 \tau _2 $. 
Let $ m $ be the line in $ \pi $ containing $ P $.  
It follows from  lemma \ref{tran12} that $ \tau_2 P \in m ,$  and thus also $ \tau _1 \tau _2 P \in m $.  
Thus $ m = l $, and it follows that $ l \in \pi $. This shows that $ \tau _1 \tau _2 \in { T }_ \pi $. $\Box$

\begin{lm}
\label{tran20}
Let $ \tau _1 $ and $ \tau _2 $ be translations with
$ \tau _1 \tau _2 \ne \tau _2 \tau _1 $.  
Then $ \tau _1 \ne 1 ,$ $ \tau _2 \ne 1 ,$ 
and $ \tau _1 $ and $ \tau _2 $ have the same direction.
\end{lm}

\noindent Proof. Choose a point $ P $ such that
$ \tau _1 \tau _2 P \ne \tau _2 \tau _1 P $. Thus $  \tau _2 P \ne \tau _1 ^ {-1} \tau _2 \tau _1 P $.  
Either $ P \ne \tau _2 P $ or $ P \ne \tau _1 ^ {-1} \tau _2 \tau _1 P $. In the first case, $ \tau _2 \ne 1 $. 
In the second case, $ \tau _1 P \ne \tau _2 \tau _1 P ,$  and again it follows that $ \tau_2 \ne 1 $.  Thus, in either case, $ \tau _2 \ne 1 $.  Similarly, $ \tau _1 \ne 1 $. 

Now set $ \tau \equiv \tau _1 \tau _2 \tau _1 ^ {-1} \tau_2 ^ {-1} $; thus $ \tau \ne 1 $. Let $ l $\, be any trace of $ \tau $. Theorem \ref{dila17} and  theorem \ref{tran18} show that $ t ( \tau _2 \tau _1 ^ {-1} \tau _2 ^ {-1} ) = t ( \tau _1^ {-1} )= t ( \tau _1 ) ,$  and it  follows from theorem \ref{tran19} that $ l \in  t ( \tau ) =  t ( \tau _1 \cdot \tau _2 \tau _1 ^ {-1} \tau _2 ^ {-1} ) = t ( \tau _1 ) $. Similarly, ~$ l  \in t ( \tau ) =  t ( \tau _1 \tau _2 \tau _1^{-1}\cdot\tau _2 ^ {-1} ) = t  ( \tau _2 ) $. Since $ \tau _1 $  and $ \tau _2 $ have the common trace $ l $, they have the same direction. $\Box$

\begin{lm}
\label{tran21}
If the translations $ \tau _1 $ 
and $ \tau _2 $ have no common trace, 
then $ \tau _1 \tau _2 = \tau _2 \tau _1 $.
\end{lm}

\begin{thm} 
\label{tran22}
Assume that for any given translation $ \tau \ne 1 ,$    there exists another translation $ \tau ' \ne 1 $ such that $ \tau $ and $ \tau ' $ have different directions.\footnote{This will follow from axiom K1 in section \ref{sec:ring}.}   Then the translation group $ T $ is commutative. 
\end{thm}

\noindent Proof. Let $ \tau _1 $ and $ \tau _2 $ be any translations, and suppose that $ \tau _1 \tau _2 \ne \tau _2 \tau _1 $.  Lemma \ref{tran20} shows that $ \tau _1 \ne 1 ,$ $ \tau _2 \ne 1 ,$ and that $ \tau _1 $ and $ \tau _2 $ have the same direction, which we will denote by $ \pi $. Thus $ \tau _1 $ and $ \tau _2 $ belong to the subgroup $ T _\pi $. Choose a translation $ \tau_3 \ne 1 $ such that $ \tau_3 $ and $ \tau _1 $ have different directions; it follows from  lemma \ref{tran21} that $ \tau_3 \tau_1 = \tau_1 \tau_3 $.

Now suppose further that $ \tau _2 \tau_3 $ and $ \tau _1 $ have a common trace; it follows that $ \tau_2 \tau_3 \in T _\pi $. Since $ \tau_3 = \tau_2 ^ {-1}\cdot\tau_2 \tau_3 ,$ we then have $ \tau_3 \in T _\pi ,$ a contradiction. Hence $ \tau_2 \tau_3 $ and $ \tau_1 $ have no common trace.

 Now, again by lemma \ref{tran21},  we have $ 
\tau_1 \cdot \tau_2 \tau_3 = \tau_2 \tau_3 \cdot \tau_1  = \tau_2 \cdot \tau_3 \tau_1 = \tau_2 \cdot \tau_1 \tau_3 $, and thus $ \tau_1 \tau_2 = \tau_2 \tau_1 ,$ a contradiction. This shows that 
$ \tau_1 \tau_2 = \tau_2 \tau_1 $. $\Box$

\section{Division ring}
\label{sec:ring}

The classical theory  of the division ring of scalars is highly nonconstructive. The main step in the  constructivization, theorem \ref{ring18}, will require both the displacement property of translations, obtained above in theorem \ref{tran9}, and the injective property of nonzero trace-preserving homomorphisms, derived below in theorem \ref{ring11}. \\

 \noindent \textbf{Axiom K1.}
\label{axK1}\textit{Given any points 
$ P $ and $ Q ,$
there exists a translation 
that maps $ P $ into $ Q $.} \\

The translation resulting from axiom K1 will be denoted $ \tau_{PQ} $. \label{taupq}

\begin{prop}
\label{ring3} \hspace*{5mm}\\
\hspace*{5mm}(a) For any points $ P $ and  
$ Q ,$ the translation
$ \tau_{PQ} $ is unique.\\
\hspace*{5mm}(b) The translation group $ T $ is commutative.
\end{prop}

\noindent Proof. The uniqueness of $ \tau_{PQ} $ follows from theorem \ref{tran17}. The  commutativity of the group $ T $ follows from  theorem \ref{tran22}. $\Box$

\begin{defn}  
\label{ring4}
\textnormal
{A map $ \alpha : T \rightarrow T $ will be called a \textit{trace-preserving homomorphism} if:\\
\hspace*{5mm}(a) For any translations 
$ \tau_1 $ and $  \tau_2 ,$ \hspace*{1mm}  
$ (\tau_1 \tau_2)^ \alpha = \tau_1^ \alpha \tau_2^ \alpha $.\\
\hspace*{5mm}(b) For any translation $ \tau ,$ \hspace*{1mm} $  t(\tau ^ \alpha) \subseteq t(\tau) $.\footnote{The inclusion in part (b) of definition 
\ref{ring4} is the reverse of that used in [A]. This is required because of the restricted notion of ``trace'' adopted in definition \ref{dila16}.}\\
The set of all trace-preserving homomorphisms will be denoted $ k $.} 
\end{defn} 

\noindent \textbf{Examples.}
\label{ring6} These examples include certain trace-preserving homomorphisms which will be required in the ring $ k $. The verifications are straightforward.

(a) The trace-preserving homomorphism denoted by $ 0 $ maps any translation $ \tau $ into the identity in $ T $. 
Thus $ \tau^0 = 1 $ for all
$ \tau \in T $.

(b) The trace-preserving homomorphism denoted by $ 1 $ is the identity map.
Thus $ \tau ^ 1 = \tau $ for all $ \tau \in T $. Since we have axiom K1, it follows that $ 1 \ne 0. $ 

(c) The trace-preserving homomorphism denoted by $ -1 $ maps any translation $ \tau $ into its inverse $ \tau^{-1} $.

(d) Let $ \sigma $ be a dilatation. The trace-preserving homomorphism denoted by
$ \alpha_ \sigma $ \label{alphasigma} is defined by 
$ \tau ^{ \alpha_ \sigma } = \sigma \tau \sigma ^ {-1} $ for all $ \tau \in T $. It is easily verified that 
$ \alpha _ \sigma \ne 0 $.

\begin{prop}
\label{ring7}
Let $ \tau $ be a translation, let $ l $ be a trace of $ \tau ,$ and let $ \alpha \in k $. If $ P \in l ,$ then also
$ \tau ^ \alpha P \in l $.
\end{prop}

\noindent Proof. This follows from 
lemma \ref{tran12}. $\Box$

\begin{defn}
\label{ring8}
\textnormal{Let $ \alpha, \beta \in k $.}
\textnormal{\\
 \hspace*{5mm}(a) Define a map, denoted $ \alpha + \beta ,$ by $ \tau ^ { \alpha + \beta } = \tau ^ \alpha \tau ^ \beta $ for all $ \tau \in T $.\\
\hspace*{5mm}(b) Define a map, denoted $ \alpha \beta ,$ by
$ \tau ^ { \alpha \beta } 
= ( \tau ^ \beta ) ^ \alpha $ for all  $ \tau \in T $. } 
\end{defn}

\begin{thm} \hspace*{5mm}\\
\label{ring9}
\hspace*{5mm}(a) For any $ \alpha, \beta  \in k ,$ the maps
$ \alpha + \beta $ and $ \alpha \beta $ are trace-preserving homomorphisms.\\
\hspace*{5mm}(b) Under definitions \ref{ring8}, 
$ k $ is a ring with identity. 
\end{thm}

\noindent Proof. The algebraic portions of the proof are straightforward; we need only examine the traces. Let $ \alpha, \beta  \in k $ and let $ \tau \in T $. 

Let $ l \in t ( \tau ^ { \alpha + \beta }) $. It follows that $ \tau ^ \alpha \tau ^ \beta \ne 1 $. Using  proposition \ref{dila15}, we have either $ \tau ^ \alpha  \ne 1 $ or $  \tau ^ \beta \ne 1 $. In the first case, $ \emptyset \ne t( \tau ^ \alpha) \subseteq 
t ( \tau ) $; thus $ \tau \ne 1 $. The second case is similar; thus we may let $ \pi $ denote the pencil $  t ( \tau ) $. It follows from theorem \ref{tran19}  that
$ \tau ^ \alpha \tau ^ \beta \in T _ \pi ,$ and thus $ l \in t ( \tau ) $. This shows that the map $  \alpha + \beta  $ has the trace-preserving property. 

Also, $ t ( \tau ^ { \alpha \beta }) = t (( \tau ^ \beta ) ^ \alpha ) \subseteq t ( \tau ^ \beta ) \subseteq t( \tau ) $; thus $ \alpha \beta \in k $. $\Box$\\

We will need to know that the product of translations with different directions is distinct from the identity. Further, we will require  this conclusion even in a situation where one of the translations is not known to be  $ \ne 1 $, only that its traces, if any, are distinct from those of the other translation. 

\begin{lm}
\label{ring10}
Let  $ \tau _1 \ne 1 $ be a translation with direction $ \pi _ 1 $, and let $ \tau  _ 2 $ be a translation with
$ t ( \tau  _ 2 ) \subseteq \pi _ 2 ,$ 
where $ \pi _ 2 $ is a pencil of lines distinct from $ \pi _ 1 $. Then 
$ \tau _ 1 \tau _ 2 \ne 1 $. 
\end{lm}

\noindent Proof. Choose a point $ P $ such that $ l _ 1 \equiv P + \tau _ 1 P $ is a trace of $ \tau _1 $. Denote by $ l _2 $ the line in $ \pi _2 $ through $ \tau _ 1 P $. It follows from lemma \ref{tran12} that $ \tau _ 2 \tau _ 1 P \in l _2 $. Since $ P \ne \tau _ 1 P = l _ 1 \cap l _ 2 ,$ it follows from axiom L1 that $ P \notin l _2 $.
Thus $ P \ne \tau _ 2 \tau _ 1 P ;$ this shows that $ \tau _ 2 \tau _ 1 \ne 1 $. $\Box$\\

Classically, the following theorem is a consequence of [A, Theorem 2.12]. The latter is proved nonconstructively using multiplicative inverses in $ k ,$ which are also derived nonconstructively. For a constructive proof that $ k $ is a division ring, we must derive theorem \ref{ring11} first, directly from the properties of translations and traces. 

\begin{thm}
\label{ring11}
Let $ \alpha $ be a trace-preserving homomorphism. If $ \alpha \ne 0 ,$ then $ \tau ^ \alpha \ne 1 $ for all translations $ \tau \ne 1 $. 
\end{thm}

\noindent Proof. Since  $ \alpha \ne 0 ,$ we may choose one translation $ \tau _ 1 $ such that $ \tau _ 1 ^ \alpha \ne 1 $. Thus $ \tau _ 1 \ne 1 $; let $ \pi _ 1 $ denote the direction of $ \tau _ 1 $. Choose any point $ P$ such that $ \tau _1 P \ne P $, and   set
$ l _ 1 \equiv P + \tau _ 1 P $; thus $ l _ 1 \in \pi _1 $. 

First consider the special case in which a translation $ \tau \ne 1 $ has direction $ \pi $ distinct from the pencil $ \pi _1 $. Using theorem \ref{tran9}, we may set $ l \equiv \tau _1 P + \tau \tau _1 P $; thus $ l \in \pi $. It follows from  lemma \ref{ring10} that $ \tau ' \equiv \tau \tau _1 \ne 1 $; let $ \pi ' $ denote the direction of $ \tau ' $. By theorem \ref{tran9} again, we may set $ l ' \equiv P + \tau ' P $; thus $ l' \in \pi ' $. Since $ \tau ' P \ne \tau _1 P = l \cap l _1 $, it follows from axiom L1 that $ \tau ' P \notin  l _1 $. Thus $ l _1 \nparallel l ' ,$ and the pencils $ \pi _1 $ and $ \pi ' $ are distinct. Since $ t (\tau _1 ^ {- \alpha }) = t ( \tau _ 1 ^ \alpha ) = \pi _1 ,$ and $ t ( \tau ^ \alpha \tau _1 ^ \alpha ) \subseteq t ( \tau \tau _1 ) = \pi ' $, it follows from lemma \ref{ring10} that $ \tau ^ \alpha \tau _1 ^ \alpha \cdot \tau _1 ^ {- \alpha } \ne 1 $; thus  $ \tau ^ \alpha \ne 1 $.

The special case shows that $ \tau ^ \alpha \ne 1 $ for any translation $ \tau \ne 1 $ with direction distinct from $ \pi _ 1 $. Use axiom K1 to construct one such translation $ \tau _2 ,$ with direction $ \pi _2 $. Thus $ \tau _2 ^ \alpha \ne 1 ,$ and it follows from the special case that $ \tau ^ \alpha \ne 1 $ for any translation $ \tau \ne 1 $ with direction distinct from $ \pi _2 $. Now consider any translation $ \tau \ne 1 $. It follows from  axiom L2 that its direction is either distinct from $ \pi _1 $ or distinct from $ \pi _ 2 $; thus $ \tau ^ \alpha \ne 1 $. $\Box$

\begin{cor}
\label{tphinj}
Let $ \alpha \ne 0 $ be a trace-preserving homomorphism. If $ \tau_1 $ and $ \tau_2 $ are translations with  $ \tau_1 \ne  \tau_2 $, then $ \tau_1^ \alpha  \ne  \tau_2^ \alpha  $. 
\end{cor}

\begin{cor}
\label{ring12}
The product of two translations with different directions has a third direction, distinct from the first two. 
\end{cor}

\noindent Proof. In the proof of the theorem, translations $ \tau _1 $ and $ \tau $ are arbitrary, with distinct directions $ \pi _ 1 $ and $ \pi $. The product $ \tau \tau _1 $ has direction $ \pi ' ,$ and the proof shows that $ \pi ' $ is distinct from $ \pi _1 $. Also, since $  P \ne \tau _1 P = l \cap l _1 ,$ it follows that $ P \notin  l $; thus $ l \nparallel l ' $. This shows that $ \pi ' $ is also distinct from $ \pi $. $\Box$

\begin{cor}
\label{ring13}
Let $ \tau _1 \ne 1 $ and $ \tau _2 \ne 1 $ be translations with different directions, and let $ \alpha $ and $ \beta $ be elements 
of $ k ,$ with $ \alpha \ne 0 $. Then 
$ \tau _1^\alpha \tau _2^\beta \ne 1 $. 
\end{cor}

\noindent Proof. It follows from theorem \ref{ring11} 
that $ \tau _1 ^ \alpha \ne 1 $; thus lemma \ref{ring10} applies. $\Box$ \\

The following  corollary is a constructive version of [A, Theorem 2.12], ``If $ \tau ^ \alpha = 1 ,$ then either $ \alpha = 0 $ or $ \tau = 1 $,'' which is constructively invalid.\footnote{See example \ref{brouF}.} The most essential constructive substitute, however, is theorem \ref{ring11} itself. 

\begin{cor}
\label{ring14} 
Let $ \alpha, \beta \in k $.

(a) Let $ \alpha \ne 0 $. If $ \tau ^ \alpha = 1 $ for some translation $ \tau ,$ then $ \tau = 1 $.

(b) If $ \tau ^ \alpha = 1 $ for some translation $ \tau \ne 1 ,$ then $ \alpha = 0 $.

(c) If $ \tau^\alpha = \tau^\beta $ for some translation $ \tau \ne 1 ,$ then $ \alpha = \beta $. Thus a trace-preserving homomorphism is uniquely determined by the image of a single translation $ \ne 1 $.

(d) If $ \tau ^ \alpha = \tau $ for some translation $ \tau \ne 1 ,$ then $ \alpha = 1 $. 

\end{cor}  

\begin{thm}
\label{ring18}
Let $ \alpha $ be a trace-preserving homomorphism with $ \alpha \ne 0 ,$ and let $ P $ be any point. Then there exists a unique dilatation $ \sigma $ with fixed point $ P $ such that  
$ \alpha = \alpha _ \sigma $. 
\end{thm}

\noindent Proof. We first prove the uniqueness, and determine a working definition for a map 
$ \sigma $. Let $ \sigma $ be as specified, and let $ Q $ be any point.
It follows from axiom K1 that
$ \sigma Q = \sigma \tau _ {PQ} P = \sigma \tau _ {PQ} \sigma ^ {-1} P = \tau _ {PQ} ^ {\alpha _\sigma} P = \tau _ {PQ} ^ \alpha P $.
This shows that $ \sigma ,$ if it exists, is unique. 

Define a map $ \sigma :\mathscr{P}  \rightarrow \mathscr{P} $ by

\begin{displaymath}
\sigma Q = \tau _ {PQ} ^ \alpha P\hspace*{10mm} \mathrm{for \; all} \hspace*{2mm} Q \in \mathscr{P} 
\end{displaymath}
Clearly, $ \sigma P = P $. 

Let $ Q $ and $ R $ be any points with $ Q \ne R $.
Then 

\begin{displaymath}
\sigma R = \tau _ {PR} ^ \alpha P = \tau _ {QR} ^ \alpha \tau _ {PQ} ^ \alpha P = \tau _ {QR} ^ \alpha \sigma Q
\end{displaymath}
Since $ \alpha \ne 0 $ and $ \tau _ {QR} \ne 1 ,$ it follows from theorem \ref{ring11}  that
$ \tau _ {QR} ^ \alpha \ne 1 $.
Theorem \ref{tran9} shows that $ \tau _ {QR} ^ \alpha \sigma Q \ne \sigma Q $;
thus $ \sigma R \ne \sigma Q $. This shows that 
 $ \sigma $ is injective.
 
Since $ Q + R $ is a trace 
of $ \tau _ {QR} ,$
~$ \sigma Q + \sigma R $ is a trace 
of $ \tau _ {QR} ^ \alpha $, and $ \alpha $ is trace-preserving, we have   
$ Q + R \parallel \sigma Q + \sigma R $. 
It now follows from 
theorem \ref{dila7} that $ \sigma $ is a dilatation.

Now let $ \tau $ be any translation, and set $ S \equiv \tau P $. Then $ \tau ^ {\alpha _ \sigma } P = \sigma \tau \sigma ^ {-1} P = \sigma \tau P = \sigma S = \tau _ {PS} ^ \alpha P =\tau ^ \alpha P $. Since  $ \tau ^ { \alpha _ \sigma  } $ and $ \tau ^ \alpha $ agree at the point $ P ,$ it follows from theorem \ref{tran17}  that these translations are equal. Thus  $ \alpha = \alpha _ \sigma $. $\Box$

\begin{defn}
\label{ring19}
\textnormal{A ring $ k $ with identity is a \textit{division ring} if for any elements $ x $ and $ y $ in $  k ,$ ~$ x  \ne y \textnormal{~~if and only if~~} x - y  \textnormal{~~is a unit in~~}  k $.}\footnote{See [MRR; \S\,II.2].} 
\end{defn}

 \noindent \textit{Note.} Let $ k $ be a ring with identity, and let $ k $ have an inequality relation that is invariant with respect to addition. Then $ k $ is a division ring if and only if, for any element $ x \in  k $, $ x \ne 0 $ if and only if $ x $ is a unit in $ k $. This applies to the ring $ k $ of trace-preserving homomorphisms. 

\begin{thm}
\label{ring20}
The ring $ k $ of trace-preserving homomorphisms has the following properties, where 
$ \alpha, ~\beta  \in k $. 

(a) If $ \alpha \ne \beta ,$ then for any $ \gamma \in k ,$ either $ \gamma \ne \alpha ,$ or  $ \gamma \ne \beta $.

(b) If $ \alpha \beta \ne 0 ,$ then $ \alpha \ne 0 $ and $ \beta \ne 0 $.

(c) $ k $ is a division ring.
\end{thm}

\noindent Proof. (a) Choose a translation $ \tau  ,$ and then a point $ P ,$ such that $ \tau ^ \alpha  P \ne \tau ^ \beta  P .$ Now, either $ \tau ^ \gamma   P \ne \tau ^ \alpha  P ,$ or $ \tau ^ \gamma  P \ne \tau ^ \beta  P ,$  and it follows that either $ \gamma \ne \alpha ,$ or  $ \gamma \ne \beta $.

(b) Choose a translation $ \tau $ such that $ \tau ^ {\alpha \beta} \ne 1 $. Then $ \emptyset \ne t(\tau ^ {\alpha \beta}) = t((\tau ^ \beta) ^ \alpha ) \subseteq t(\tau ^ \beta) $. It follows that $ (\tau ^ \beta) ^ \alpha \ne 1 $ and $ \tau ^ \beta \ne 1 $; thus $ \alpha \ne 0 $ and $ \beta \ne 0 $. 

(c) Let $ \alpha \in k $ with $ \alpha \ne 0 $. Using theorem \ref{ring18}, construct a dilatation $ \sigma $ such that $ \alpha = \alpha _ \sigma $. It is  clear that $ \alpha _ {\sigma ^ {-1}} \alpha _ \sigma = \alpha _ \sigma \alpha _ {\sigma ^ {-1}} = 1 $; thus $ \alpha _ {\sigma ^ {-1}} $ is the inverse of  $ \alpha  $. Conversely, let $ \alpha $ be a unit in $ k $; then $ \alpha \beta = 1 $ for some  $ \beta \in k $. It then follows from condition (b) that 
$ \alpha \ne 0 $. $\Box$

\section{Coordinates}
\label{sec:coor}

Whereas axiom K1 provided  translations mapping any point to any point, axiom K2 will provide  dilatations that expand about a fixed central point, mapping other points arbitrarily along radials. These dilatations will then lead to the required scalars and coordinates.\\

 \noindent \textbf{Axiom K2.}\label{axK2} \textit{Let $ P $ be any point. If ~$ Q $ and $ R $ are points collinear with $ P,$ and each is distinct from $ P,$ then there exists a dilatation $ \sigma $ with fixed point $ P $ that maps $ Q $ into $ R $.}\\

It follows from theorem \ref{dila14} that the dilatation $ \sigma $ resulting from axiom K2 is unique. 

\begin{thm}
\label{coor3} 
The following are equivalent.

(a) Axiom K2.

(b) If $ \tau _1 \ne 1 $ and $ \tau _2 \ne 1 $ are translations with the same direction, then there exists a unique trace-preserving homomorphism $ \alpha \ne 0 $ in $ k $ such that $ \tau _2 = \tau _1 ^  \alpha $.
\end{thm}

\noindent Proof. Given axiom K2, let $ \tau _1 $ and 
$ \tau _2 $ satisfy the hypotheses in (b). Choose any point $ P ,$  set $ Q \equiv \tau _1 P $,  and set 
$ R \equiv \tau _2 P $. Use the axiom  to construct a dilatation $ \sigma $ with fixed point $ P $ such that
$ \sigma Q = R ,$ and set $ \alpha \equiv \alpha _ \sigma $. Since the translations $ \tau _2 $ and 
$ \sigma \tau _1 \sigma ^ {-1} $ agree at the point $ P ,$ it follows from theorem \ref{tran17}  that they are equal; thus $ \tau _2 = \tau _1 ^ \alpha $. The uniqueness of $ \alpha $ follows from proposition \ref{ring14}(c).

Given (b), let $ P $ be any point, and let the points $ Q $  and $ R $ satisfy the hypotheses in axiom K2. Use axiom K1 to construct the translations $ \tau _ {PQ} $ and $ \tau _ {PR} $; thus $ \tau _ {PQ} \ne 1 ,$ $ \tau _ {PR} \ne 1 ,$ and these translations have the same direction. Let $ \alpha \ne 0 $ be the element of $ k $ such that $ \tau _ {PR} = \tau _ {PQ} ^ \alpha $. Use theorem \ref{ring18}  to construct the dilatation $ \sigma $ with fixed point $ P $ such that $ \alpha = \alpha _ \sigma $. It then follows that $ R = \tau _ {PR} P = \sigma \tau _ {PQ} \sigma^{-1}P =  \sigma Q $. $ \Box $  

\begin{prop}
\label{coor3A}
If the statement  of axiom K2 holds at a single point $ P $, then it holds at every point, and thus axiom K2 is valid. 
\end{prop}

\noindent Proof. Let the statement  hold at the 
point $ S ,$ and let the points $ P, \, Q $ and $ R $ be as in axiom K2. Using axiom K1, construct the translation $ \tau \equiv \tau _{PS} $. Then  $ P + Q \parallel  S + \tau Q $, and $ P + R \parallel  S + \tau R $. It follows that  $ \tau Q $ and $ \tau R $ are collinear with $ S ,$ and distinct from $ S $. Let $ \sigma _1 $ be the dilatation with fixed point $ S $ such that $ \sigma _1 \tau Q = \tau R .$ It then follows that $ \sigma \equiv  \tau ^{-1} \sigma _1 \tau $ is the required dilatation. $\Box$

\begin{thm}
\label{coor4}
Let $ \tau _1 \ne 1 $ be a translation. For any translation $ \tau _ 2 $ with $ t ( \tau _ 2) \subseteq t ( \tau _ 1) ,$ there exists a unique element $ \alpha $ in $ k $ such that $ \tau _2 = \tau _1 ^ \alpha $. 
\end{thm}

\noindent Proof. Choose any point $ P $ such that  $ \tau _1 P \ne P $. It then follows that either
$ \tau _2 P \ne P $ or $ \tau _2 P \ne \tau _ 1 P $. In the first case, $ \tau _2 \ne 1 ,$ and theorem \ref{coor3} applies directly. In the second case, $ \tau _2 \ne \tau _ 1 $ and it follows from theorem \ref{tran19} that $ \tau _1 ^ {-1} \tau _2 $ has the same direction as $ \tau _1 $. Use theorem \ref{coor3} to construct an element $ \beta $ in $ k $ such that $ \tau _1 ^ {-1} \tau _2 = \tau _ 1 ^ \beta $; thus $ \tau _2 = \tau _ 1^ {\beta+1} $. The uniqueness follows from corollary  \ref{ring14}(c). $\Box$\\ 

Theorem \ref{coor4} is a constructive substitute for  [A; page 63, Remark]. This remark, an essential part of the classical theory, requires the nonconstructive statement: ``either $ \tau _2 = 1 $ or $ \tau _2 \ne 1 $''.\footnote{See example \ref{brouE}.} The theorem here covers all cases, without determining whether or not $ \tau _2 = 1 $. 

\begin{thm}
\label{coor5}
Let $ \tau _1 \ne 1 $ and $ \tau _2 \ne 1 $ be translations with different directions. For any 
translation $ \tau ,$ there exist unique 
elements $ \alpha $ and $ \beta $ in $ k $ such that

\begin{displaymath}
\tau = \tau _1 ^ \alpha \tau _2 ^ \beta
\end{displaymath}
If $ \tau \ne 1 ,$ then either $ \alpha \ne 0 $ or $ \beta \ne 0 ,$ and conversely.
\end{thm}

\noindent Proof. Choose any point $ P ,$ and 
set $ Q \equiv \tau P $. Let $ l _2 $ be the $ \tau _2 $ trace through $ P ,$  let $ l _ 1 $ be the $ \tau _1 $ trace through $ Q $, and set 
$ R \equiv l _1 \cap l _2 $. Let $ l $ be a trace of $ \tau _ {PR} $. Then $ P \ne R $, and $ l _2 $ is also a trace of $ \tau _ {PR} $. Thus $ l \parallel l_2 $, and this shows that $ t ( \tau _ {PR} ) \subseteq t ( \tau _2 ). $ Similarly,  $ t ( \tau _ {RQ} ) \subseteq t ( \tau _1 ). $ Using theorem \ref{coor4}, construct elements $ \alpha $ and $ \beta $ in $ k $ such that $ \tau _ {PR} = \tau _2 ^ \beta $ and $ \tau _ {RQ} = \tau _1 ^ \alpha $. It follows that $  \tau _1 ^ \alpha \tau _2 ^ \beta $ takes $ P $ into $ Q $, and it is therefore equal to $ \tau $.

Now let $ \tau _1 ^ \alpha \tau _2 ^ \beta = \tau _1 ^ \gamma  \tau _2 ^ \delta ,$ and set $ \tau _ 3 \equiv \tau _1 ^ { \alpha - \gamma} = \tau _2 ^ { \delta - \beta } $. Since $ t ( \tau _ 3) \subseteq 
t ( \tau _ 1) \cap t ( \tau _ 2) = \emptyset ,$ it follows that $ \tau _ 3 = 1 ,$ and it then follows from corollary  
\ref{ring14}(b)  that $ \alpha - \gamma = \delta - \beta = 0 $. Thus $ \alpha $ and $ \beta $ are unique.

Finally, let $ \tau \ne 1 $. By proposition 
\ref{dila15}, either $ \tau _1 ^ \alpha \ne 1 $ 
or $ \tau _ 2 ^ \beta \ne 1 $. Thus either $ \alpha \ne 0 $ or $ \beta \ne 0 $. The converse follows from corollary \ref{ring13}. $\Box$\\

We are now prepared to place the capstone to the coordinatization theory. 

\begin{thm}
\label{coor6}
Let $ \mathscr{G} = ( \mathscr{P}, \mathscr{L} ) $ be a Desarguesian plane. Select a point $ O ,$ and select translations $ \tau _1 \ne 1 $ and $  \tau _2 \ne 1 $ with different directions. 

(a) To any point $ P $ there corresponds a unique coordinate pair $ (x,y) $ in $ k ^ 2 $ such that
\begin{displaymath}
\tau _ {OP} = \tau _1 ^ x \tau _2 ^ y
\end{displaymath}

(b) The resulting map $ \mathscr{P} \rightarrow k ^ 2 $ is a bijection. Thus the  inequality relations on $ \mathscr{P} $ and  on $ k^2 $  correspond under this map. 

(c) To any line $ l ,$ there correspond elements 
$ \alpha, \beta, \gamma, \delta $ in $ k ,$ with either 
$ \gamma \ne 0 $ or $ \delta \ne 0 ,$ such that the points on $ l $ are the points with coordinates in the set
\begin{displaymath}
L = \{(\alpha + t\gamma, \beta + t\delta ) : 
t \in k \}
\end{displaymath}
Conversely, if  elements $ \alpha, \beta, \gamma, \delta $ in $ k $ are given, with either $ \gamma \ne 0 $ or $ \delta \ne 0 ,$ then the set of points with coordinates in the set $ L $  determines a line in $ \mathscr{L} $. 

(d) The principal relation $ P \notin l $ of definition \ref{ax-2b} corresponds to the condition ``$ (x,y) \ne 
(\alpha + t\gamma, \beta + t\delta ) $ for all $ t \in k $.''  
\end{thm}

\noindent Proof. (a) This follows from theorem \ref{coor5}. 

(b) For any pair $ (x,y) \in k ^ 2 ,$ set $ P \equiv \tau _1 ^ x \tau _2 ^ y O $; thus $ P \rightarrow (x,y) $. This shows that  the map is onto $ k ^ 2 $. Now let $ P $ and $ Q $ be points with $ P \rightarrow (x,y) $ and $ Q \rightarrow (z,w) $. If $ P $ and $ Q $ are distinct, then $ \tau _ {PQ} \ne 1 $ and $ \tau _ {PQ} = \tau _ {OQ} \tau _{PO} = \tau _1 ^ z \tau _2 ^ w \tau _1 ^ {-x} \tau _2 ^ {-y} = \tau _1 ^ {z-x} \tau _2 ^ {w-y} $. It follows from proposition \ref{dila15} that either $ \tau _1 ^ {z-x} \ne 1 $ or $ \tau _2 ^ {w-y} \ne 1 $; thus either $ z \ne x $ or $ w \ne y ,$ and hence $ (x,y) \ne (z,w) $. This shows that the map $ \mathscr{P} \rightarrow k ^ 2 $ is injective. Conversely, if $ (x,y) \ne (z,w) ,$ then corollary \ref{ring13} and a reversal of the last argument will show that  $ P \ne Q $. Thus the inverse map is injective. 

(c) Given a line $ l ,$ choose any point $ P $ on $ l ,$ with coordinates $ ( \alpha, \beta ) ;$ thus  $ \tau _ {OP} = \tau _1 ^ \alpha  \tau _2 ^ \beta  $. Choose any translation $ \tau $ with $ l $ as a trace. Use 
theorem \ref{coor5} to construct elements $  \gamma $ and $ \delta $ in  $ k $ such that $ \tau = \tau _1 ^ \gamma \tau _2 ^ \delta $; thus either 
$ \gamma \ne 0 $ or $ \delta \ne 0 $. For any point $ Q $ on $ l ,$ use theorem \ref{coor4} to construct the unique element $ t \in k $ such that $ \tau _{PQ} = \tau ^ t $. Then $ \tau _{OQ} = \tau _{PQ}\tau _{OP} = \tau^t \tau _{OP} = \tau _ 1 ^ {t \gamma + \alpha } \tau _ 2 ^ {t \delta + \beta } $, and thus $ Q $ has coordinates in the set $ L $. Conversely, a reversal of this argument will show that if a point has  coordinates in the set $ L $, then it lies on $ l $. 

Now let $ \alpha, \beta, \gamma, \delta $ be elements of $ k ,$ with either $ \gamma \ne 0 $ or $ \delta \ne 0 $. Let  $ P $ be the point with coordinates $ ( \alpha, \beta) $, and set $ \tau \equiv \tau _ 1 ^ \gamma \tau _2 ^ \delta $. It follows from corollary \ref{ring13} that $ \tau \ne 1 $; let $ l $ denote the trace of $ \tau $ containing $ P $. The construction above shows that the points on $ l $ are those with coordinates in the set L. 

(d) This follows from part (b). $ \Box $\\ 

Using the expression for the set $ L $ in theorem \ref{coor6}, one may obtain parametric equations for a line $ l ,$ and, if $ k $ is commutative, also an equation in the  form $ ax + by + c = 0 $, where either $ a \ne 0 $ or $ b \ne 0 $.

\section{Desargues}
\label{sec:des} 

Assuming now only the axioms in groups \textbf{G}  and \textbf{L}, we demonstrate that the axioms in group \textbf{K} are equivalent to \textit{Desargues's Theorem}; this theorem has two variations, stated below as postulates D1 and D2. Using Desargues's theorem as an alternative to axiom group \textbf{K} would have the advantage that these postulates involve only direct properties of the parallelism concept.\\

 \noindent \textbf{Postulate D1.}\label{postD1} 
\textit{Let $ l _1, ~l_2, ~l_3 $ be distinct parallel lines. Let $ P,~P' \in l_1 $; ~$ Q,~Q' \in l_2 $; ~and ~$ R,~R' \in l_3 $. If} 

\begin{displaymath}
 P + Q \parallel P' + Q'  
\textnormal{~~and~~}  P + R \parallel P' + R' 
\end{displaymath}
\textit{then}  
\begin{displaymath}
Q + R \parallel Q' + R'
\end{displaymath}

\vspace*{4mm}

 \noindent \textbf{Postulate D2.}\label{postD2} 
\textit{Let $ l _1, ~l_2, ~l_3 $ be distinct concurrent lines. Let 
$ P,~P' \in l_1 $; ~$ Q,~Q' \in l_2 $; 
and ~$ R,~R' \in l_3 $; with these points each distinct from the point of concurrence.
If~ $ P + Q \parallel P' + Q'  ~and~ P + R 
\parallel P' + R' ,$ ~then~ 
$ Q + R \parallel Q' + R' $.}

\begin{thm}
\label{des4}
Axiom K1 implies postulate D1. Axiom K2 implies postulate D2.
\end{thm}

\noindent Proof. Assume axiom K1, consider the configuration of postulate D1, and set $  \tau \equiv \tau_ {PP'} $. Then $ P' + Q' \parallel P + Q \parallel \tau P + \tau Q = P' + \tau Q, $ and thus $ \tau Q \in P' + Q' $. The traces of the translation $ \tau ,$ if any, are contained in the pencil $ \pi $ that contains the three lines $ l_i $. Thus lemma \ref{tran12}  applies, and $ \tau Q \in l_2 $. It follows that $ \tau Q = l_2 ~\cap~ P' + Q' = Q' $. Similarly, $ \tau R = R' $. Thus $Q + R \parallel \tau Q + \tau R = Q' + R' $. 

Now assume axiom K2 and consider the configuration of postulate D2. The proof is similar to the proof of D1. Axiom K2 provides a dilatation that has the common point $ V $ of the given lines as fixed point, and  that maps $ P $ into $ P' $. Lemma \ref{dila23} is used now in lieu of lemma \ref{tran12}. $\Box$

\begin{lm}
\label{des5}
Let $ P $ and $ P' $ be distinct points, and set $ l \equiv P + P' $. For any point $ Q $ outside $ l ,$ 
let $ l' $ denote the line through $ Q $ that is parallel to $ l $, set $ m \equiv P + Q ,$ and 
let $ m' $ denote the line through $ P' $ that is parallel to $ m $. Then $ l' \nparallel m' $. If we set $ Q' \equiv l' \cap m' ,$ then $ Q' \ne Q, ~~Q' \ne P', $ and ~$ P + Q \parallel P' + Q' $. 
\end{lm}

\noindent Proof. Since $ Q \notin l, $ we have $ l \nparallel m; $ thus $ l' \nparallel m'. $ Since $ P' \ne P = l \cap m, $ it follows from axiom L1 that $ P' \notin m. $ Thus $ m \ne m', $ and it follows from theorem \ref{ax35}  that $ Q' \ne Q. $ Since also $ l \ne l', $ we have $ Q' \ne P'. $ Finally, $ P + Q = m \parallel m' = P' + Q' $. $\Box$

\begin{defn}
\label{des6}
\textnormal{The map $ Q \rightarrow Q' $ constructed in lemma \ref{des5}  will be called a \textit{partial translation}, and will be denoted $ \lambda _ {PP'} $. We extend the definition of ``trace'' to these maps.}
\end{defn}

\begin{lm}
\label{des7}
Let $ P $ and $ P' $ be distinct points, and set $ l \equiv P + P' $.
Consider the map defined by
\begin{displaymath}
Q' \equiv \lambda _ {PP'} Q \qquad
for ~all~ Q \notin l
\end{displaymath}

(a) The map $ \lambda _ {PP'} $ is defined at all points $ Q $ outside $ l $. In this domain, the map is injective.

(b) The traces of the map $ \lambda _ {PP'} $ are parallel to $ l $.

(c) For any point $ Q $ outside 
$ l ,$ ~$ \lambda _ {QQ'} P = P' $.

(d) Let $ Q $ and $ R $ be distinct points outside $ l $. If postulate D1 is valid, then $ Q + R \parallel Q' + R' $.

(e) Let $ Q $ be a point outside $ l $. If postulate D1 is valid, then the maps $ \lambda _ {PP'} $ and $ \lambda _ {QQ'} $ agree at all points in their common domain. 
\end{lm}

\noindent Proof. (a) Let $ Q $ and $ R $ be distinct points outside $ l $. Thus $ Q' $ is determined by $ l',~ m $ and $ m',$ as constructed  in lemma \ref{des5}. Similarly, let $ R' $ be determined 
by $ l'',~ n, $ and $ n' $.

Since $ R \ne Q = l' \cap m, $ it follows from axiom L1 that either $ R \notin l' $ or $ R \notin m $. In the first case, we have $ l'' \ne l' ,$ and it follows from theorem \ref{ax35}  
that $ R' \ne Q' $.
In the second case, we have $ m \ne n, $ and thus $ m \nparallel n; $ it follows that $ m' \nparallel n'. $ Since $ R' \ne P' = m' \cap n', $ it follows that $ R' \notin m', $ and thus $ R' \ne Q' $.

(b) and (c) These are clear from the construction in lemma \ref{des5}.

(d) With the notation as in the proof of part (a), first consider the special case in which $ R \notin l'. $ Then $ l' \ne  l'', $ and we have a Desargues configuration.
Thus ~$ Q + R \parallel Q' + R' $.

In the general case, suppose that  $ Q + R \nparallel Q' + R'. $ Then the condition $ R \notin l' $ would lead, using the special case, to a contradiction; hence $ R \in l' ,$ and it follows from part (b) that also $ R' \in l' $. Thus $ Q + R = Q' + R' $, a contradiction. 
Hence $ Q + R \parallel Q' + R' $. 

(e) This leads to a Desargues configuration. $\Box$

\begin{thm}
\label{des9}
Postulate D1 implies axiom K1.
\end{thm}

\noindent Proof. (1) Assume postulate D1, and let $ P $ and $ P' $ be any points; we must construct a translation that maps $ P $ into $ P' $.

(2) Consider first the special case in which $ P $ and $ P' $ are distinct, and 
set $ l_1 \equiv P + P' $. Set $ \lambda _1 \equiv \lambda _ {PP'} $; thus $ \lambda_1 $ is defined outside $ l_1 $. Using lemma \ref{ax37}, choose one point $ P_2 $ outside $ l_1 ,$ ~set $ P_2' \equiv \lambda _1 P _2, $  ~set $ l_2 \equiv P _2 + P_2', $ and set~ $ \lambda _2 \equiv \lambda_{P_2 P_2'} $. Lemma \ref{des7} shows that $ \lambda _1 $ and $ \lambda _2 $ agree in their common domain, and that $ l_1 $ and $ l_2 $ are parallel and distinct. Define a map $ \tau $ by $ \tau Q \equiv Q' \equiv \lambda _i Q, $ whenever $ Q \notin l _i $. Lemma \ref{tran7}  shows that $ \tau $ is defined at all points. Lemma \ref{des7}  shows that $ \tau P = P', $ and that the lines $ l_Q \equiv Q + Q', $ for all points $ Q ,$ form a pencil of lines $ \pi $. 

(2.1) \textit{The map $ \tau $ is injective.}   
Let $ Q $ and $ R $ be distinct points. We may assume that $ Q \notin l_1; $ thus $ l_Q \ne l_1. $ Now, either $ R \notin l_Q, $ or $ R \notin l_1. $ In the first case, $ l_R \ne l_Q, $ and it follows from theorem \ref{ax35}  that 
$ R' \ne Q' $. In the second case, 
lemma \ref{des7}(a) applies, and again $ R' \ne Q' $.

(2.2) \textit{The map $ \tau $ is a dilatation.} To verify the direction-preserving  property, let $ Q $ and $ R $ be distinct points.

(2.2a) First consider the special case in which\\

\hspace*{15mm} (a) \textit{At least three points lie on any given line.}\\

 \noindent It then follows from theorem \ref{ax41}  that each pencil of lines contains at least three lines. Now we may construct a line $ l_3 \equiv P_3 + P_3' $ in $ \pi ,$ distinct from both $ l_1 $ and $ l_2 $. Set $ \lambda _3 \equiv \lambda _ {P_3P_3'} $; it follows that $ \tau S = \lambda _3S $ for all points $ S $ outside $ l_3 $.

Applying lemma \ref{tran7}, we find that the point $ Q $ lies outside at least two of the three lines $ l_i ,$ and that the points $ Q $ and $ R $ together lie outside one of these lines. Now lemma \ref{des7}(d) applies, and $ Q + R \parallel Q' + R' $.

(2.2b) Now consider the general case, and suppose that\\

\hspace*{15mm} (b)~ $ Q + R \nparallel Q' + R' $\\

 \noindent We may assume that $ l_Q \ne l_1$, and thus either $ R \notin l_1 $ or $ R \notin l_Q $.

In the first case, $ l_R \ne l_1 $. Suppose further that  $ l_R \ne l_Q; $ then the pencil $ \pi $ contains three distinct lines, and the special case (a) applies, a contradiction. Thus $ l_R = l_Q, $ contradicting condition (b). 

In the second case,  $ l_Q \nparallel Q + R. $ Making use of condition (b), set $ T \equiv Q + R ~\cap~ Q' + R'. $ Since $ Q \ne Q', $ either $ T \ne Q $ or $ T \ne Q'. $ In the first subcase, $ T \ne l_Q ~\cap~ Q + R; $ thus $ T \notin l_Q, $ and it follows that $ T \ne Q'. $ Thus, in either subcase, $ T \ne Q' $. Similarly, 
$ T \ne R'. $ Thus there exists a third point $ T $ on the line $ Q' + R' $, and the special case (a) applies, again contradicting condition (b).

Since condition (b) leads to a contradiction in each case, we have $ Q + R \parallel Q' + R' $. It  now follows from theorem \ref{dila7} that $ \tau $ is a dilatation.

(2.3) \textit{The map $ \tau $ is a translation.} This follows from lemma \ref{des7}(b).

(3) Now consider the general case. Choose any point $ Q $ distinct from $ P ,$ and use the special case (2) to construct the translation $ \tau _ {PQ} $. Either 
$ P' \ne P $ or $ P' \ne Q $. In the second case, the translation $ \tau \equiv \tau _{QP'} \tau _{PQ} $ takes the point $ P $ into $ P' $. This establishes axiom K1. $\Box$\\   

 \noindent \textit{Problems.} In the classical proof that postulate D1 implies axiom K1 [A, Theorem 2.17], two disjoint cases are considered: the four-point geometry, and all other geometries. For theorem \ref{des9}, we have been unable to make such a clear distinction  constructively. The proof is carried out first for the special case (a); the general case then discovers, under the assumption (b), a third point on a line. This raises the question of whether the four-point, nine-point, and larger geometries can be distinguished constructively, and the corresponding question for fields. 

\begin{lm}
\label{des25} 
Let $ l $ be a line, and let $ V $ a point on $ l $. Let $ P $ and $ P' $ be points on $ l ,$ each distinct from $ V $. For any point $ Q $ outside $ l ,$ set $ l' \equiv  V + Q ,$ set $ m \equiv P + Q ,$ and let ~$ m' $ denote the line through $ P' $ that is parallel to $ m $. 
Then $ l' \nparallel m' $. If we set $ Q' \equiv l' \cap m' ,$ then $ Q' \ne V ,$ ~$ Q' \ne P' ,$ and ~$ P + Q \parallel P' + Q' $. 
\end{lm}

\noindent Proof. Since $ Q \notin l, $ we have $ l \nparallel l'; $ it follows from axiom  L1 that $ P \notin l' .$ 
Thus $ m \nparallel l', $ and $ l' \nparallel m'. $ Since $ P' \ne V = l \cap l', $ we have $ P' \notin l' $; thus  $ P' \ne Q'. $ Since $ Q \notin l, $ we have 
$ m \nparallel l $; thus $ m' \nparallel l $. Now 
$  Q' \ne P' = l\cap m'; $ it follows that 
$ Q' \notin l $, and $ Q' \ne V. $ 
Finally, $ P + Q = m \parallel m' = P' + Q'. $ $\Box$

\begin{defn}
\label{des26}
\textnormal{The map $ Q \rightarrow Q' $ constructed in lemma \ref{des25} will be called a \textit{partial dilatation}, and will be denoted $ \lambda  _ {VPP'} $. We extend the definition of ``trace'' to these maps.}
\end{defn}

\begin{lm}
\label{des27}
Let $ l $ be a line, and let $ ~V $ a point on $ l $.
Let $ P $ and $ P' $ be points on $ l ,$ each distinct from $ V $. 
Consider the map defined by
\begin{displaymath}
Q' \equiv \lambda  _ {VPP'} Q \qquad
for ~all~ Q \notin l
\end{displaymath}

(a) The map $ \lambda  _ {VPP'} $ is defined at all points $ Q $ outside $ l $. In this domain, the map is injective.

(b) The traces of the 
map $ \lambda  _ {VPP'} $ all pass through $ V $.

(c) For any point $ Q $ outside 
$ l ,$ ~$ \lambda _ {VQQ'} P = P' $.

(d) Let $ Q $ and $ R $ be distinct points outside $ l $. If postulate D2 is valid, then $ Q + R \parallel Q' + R' $.

(e) Let $ Q $ be a point outside $ l $. If postulate D2 is valid, then the maps $ \lambda _ {VPP'} $ and $ \lambda _ {VQQ'} $ agree at all points in their common domain. 
\end{lm}

\noindent Proof. This is similar to the proof of lemma \ref{des7}. $\Box$\\

The proof that postulate D2 implies axiom K2, while similar to the proof of theorem \ref{des9}, will include several differences. The main difference is in defining a map from a collection of partial maps. In lieu of lemma \ref{tran7}, used in step (2) of theorem \ref{des9} to show that the map is everywhere defined, we must now use axiom L1, which applies only to points distinct from the point of concurrence. Thus we will require  the following extension theorem:

\begin{thm}
\label{des201}
Let $ V $ be any point, and let 
\begin{displaymath}
\mathscr{P}_V \equiv \{ Q \in \mathscr{P} : Q \ne V \} 
\end{displaymath}
be the plane punctured at $ V.$ Let 
$ \sigma _0 :\mathscr{P}_V \rightarrow \mathscr{P}_V $ be a map that is injective, has the direction-preserving  property, and has its traces all passing through $ V. $ Then $ \sigma _0 $ may be extended to a 
dilatation $ \sigma $ with fixed point $ V.$
\end{thm}

\noindent Proof. (1) Choose a point $ U $ distinct from $ V .$ Using theorem \ref{ax20} and axiom L1, construct   nonparallel 
lines $ l _1 $ and $ l_2 $ through $ U ,$ such that 
 $ V $ lies outside each line. 

(2) Let $ Q $ be any point. Either $ Q \ne V $ or $ Q \ne U $. In the first case, set $ \sigma Q \equiv \sigma _0Q. $  

(3) In the case $ Q \ne U, $ we may assume that $ Q \notin l_1. $ Choose distinct points $ P _1 $ and $ P_2 $ on $ l_1. $ Set $ l \equiv V + P _1, $ set $  l' \equiv V + P _2,  $ set $ m \equiv P _1 + Q, $ and set $ n \equiv P _2 + Q. $ Suppose $ \sigma _0 P _1 \notin l $; then $ \sigma _0 P _1 \ne P_1 $. Thus $ q \equiv P_1 + \sigma _0 P _1 $ is a trace of $ \sigma _0 $ and passes through  the point $ V $. It follows that $ q = l $, a contradiction. Hence $ \sigma _0 P _1 \in l $. Similarly, $ \sigma _0 P _2 \in l' $. Denote by $ m' $ and $ n' $ the lines parallel to $ m $ and $ n ,$ through $ \sigma _0 P _1 $ and $ \sigma _0 P _2. $ Since $ m \nparallel n, $  we have  $ m' \nparallel n'; $ set $ \sigma Q \equiv m' \cap n'. $ Thus $ \sigma $ 
is defined at all points of the plane.

(3a) If $ Q = V, $ then $ m = l, ~n = l', ~m' = l, ~n' = l', $ and $ \sigma Q = l' \cap l $. Thus $ \sigma Q = V $, and $ V $  is a fixed point of $ \sigma $. 

(3b) If $ Q \ne V, $ then $ m = P _1 + Q \parallel \sigma _0 P _1 + \sigma _0 Q $ and it follows that $ m' = \sigma _0 P _1 + \sigma _0 Q. $ Similarly, $ n' = \sigma _0 P _2 + \sigma _0 Q $; hence $ \sigma Q = \sigma _0 Q $. This shows that $ \sigma $ extends the map $ \sigma _0 $.

(4) \textit{The map $ \sigma $ is single-valued.} Let $ \sigma ' $ be a map defined by the above method, although with different choices of $ U, \,l _i, $ and $ P _i, $ and let $ Q $ be any point. Suppose that $ \sigma' Q \ne \sigma Q. $ It follows from step (3b) that $ Q = V $; by step (3a) this is a contradiction. This shows that $ \sigma' Q = \sigma Q. $

(5) \textit {If $ Q $ is any point with $ \sigma Q \ne V, $ then $ Q \ne V. $ } Either $ Q \ne V $ or $ Q \ne U. $ In the second case, we have $ \sigma Q \ne  V = l \cap l' $; thus we may assume that $ \sigma Q \notin l $. It follows that $ m' \nparallel l, $ and $ m \nparallel l. $ Since $ Q \ne P_1 = m \cap l, $ it follows that $ Q \notin l, $ and thus $ Q \ne V. $ 

(6) \textit{The map $ \sigma  $ is injective.} Let $ Q $ and $ R $ be any distinct points. Either $ V \ne Q $ or $ V \ne R; $ let us assume the latter. In this case, $ \sigma R = \sigma _0 R  \ne V $. Now, either 
$ \sigma Q \ne \sigma R ,$ or $ \sigma Q \ne V $. 
In the second subcase, step (5) shows that $  Q \ne V ,$ and hence $ \sigma Q = \sigma _0 Q \ne  \sigma _0 R = \sigma R $. 

(7) \textit{The map $ \sigma  $ has the direction-preserving  property.} Let $ Q $ and $ R $ be any distinct points, and suppose that

\begin{displaymath}
(a)~~  Q + R ~ \nparallel ~ \sigma Q + \sigma R
\end{displaymath}
Either $ V \ne Q $ or $ V \ne R; $ it suffices to consider the first case. The condition $ R \ne V $ would then imply that $ \sigma Q + \sigma R = \sigma _0 Q + \sigma _0 R, $ contradicting (a); hence $ R = V, $ and thus also $ \sigma R = V. $ Now condition (a) yields

\begin{displaymath}
(b)~~ Q + V ~ \nparallel ~ \sigma Q + V
\end{displaymath}
with $ Q $ distinct from the point of intersection $ V $. Hence $ Q \notin \sigma Q + V, $ and $ Q \ne \sigma Q. $ Using the hypothesis on traces, we have then $ V + Q = V + \sigma Q, $ contradicting (b). This shows that $ Q + R \parallel \sigma Q + \sigma R. $  

(8) It now follows from theorem \ref{dila7} that $ \sigma $ is a dilatation. $\Box$  

\begin{thm}
\label{des29}
Postulate D2 implies axiom K2. 
\end{thm}

\noindent Proof. Assume postulate D2, let $ V $ be any point, and let $ P $ and $ P' $ be points collinear with $ V $ and distinct from $ V; $ we must construct a dilatation with fixed point $ V $ that maps $ P $ into $ P'. $ 

Set $ l _1 \equiv V + P, $ and $ \lambda _1 \equiv \lambda _ {VPP'}. $ Select a point $ P_2 $ outside $ l _1, $ set $ P _2' \equiv \lambda _1 P _2, $ set $ l _2 \equiv V + P _2, $ and set $ \lambda _2 \equiv \lambda _ {VP_2 P_2'}. $ Lemma \ref{des27}(e) shows that $ \lambda _1 $ and $ \lambda _2 $ agree in their common domain. 

Define a map $ \sigma _0 : \mathscr{P} _V \rightarrow \mathscr{P} _V $ by $ \sigma _0 Q \equiv Q' \equiv \lambda _i Q $ whenever $ Q \notin l _i. $ It follows from axiom L1 that $ \sigma _0 $ is defined on $ \mathscr{P} _V. $ Lemma \ref{des27} shows that $ \sigma _0 P = P', $ and that the traces of $ \sigma _0, $ if any, pass through $ V. $ 

\textit{The map $ \sigma _0 $ is injective.}  Let $ Q $ and $ R $ be distinct points in $ \mathscr{P} _V $; set $ l _Q \equiv V + Q $ and $ l _R \equiv V + R $. We may assume that $ Q \notin l _1; $ thus $ l _Q \nparallel l_1. $ Now, either $ R \notin l _Q, $ or $ R \notin l _1. $ In the first case, $ l _R \nparallel l _Q, $ and it follows that $ R' \ne Q'. $ In the second case, lemma \ref{des27}(a) applies. 

It follows from theorem \ref{ax20} that there are at least three distinct lines through $ V. $ Choose a line $ l _3 $ through $ V, $ distinct from both $ l _1 $ and $ l _2, $ choose a point $ P _3 \ne V $ on $ l_3, $ set $ P _3' \equiv \lambda _1 P _3, $ and set $ \lambda _3 \equiv \lambda _ {VP_3P_3'}. $ It follows that $ \sigma _0S = \lambda  _3S $ whenever $ S \notin l _3. $ 

\textit{The map $ \sigma _0 $ has the direction-preserving  property.} Let $ Q $ and $ R $ be distinct points  in $ \mathscr{P} _V $. The point $ Q $ will lie outside at least two of the three lines $ l _i, $ and points $ Q $ and $ R $ together will lie outside one of these lines. Now lemma \ref{des27}(d) applies.

The map $ \sigma _0 $  now  satisfies the conditions of the extension theorem \ref{des201}; this yields the required dilatation $ \sigma, $ and establishes axiom K2. $\Box$

\section{Pappus}
\label{sec:pap} 

 Commutativity of the division ring $ k $ of trace-preserving homomorphisms will be shown equivalent to Pappus's Theorem, stated below as postulate P.\\

 \noindent \textbf{Postulate  P.}\label{pap2}  \textit{Let $ l $ and $ m $ be nonparallel lines with common point $ P $. Let $ Q, ~Q', ~Q'' $ be  points on $ l $, and let $ R, ~R', ~R'' $ be  points on $ m $, with each of these points distinct from the common point $ P $. If $ Q + R' \parallel Q' + R'' $ and $ Q' + R \parallel Q'' + R' $, then $ Q + R \parallel Q'' + R'' $.}

\begin{prop} 
\label{pap3}
Let $ P $ be any point. The following are isomorphic:

(a) The subgroup $ D _P $ of dilatations with fixed point $ P $.

(b) The multiplicative group $ k^* $ of non-zero elements in the division ring $ k $ of trace-preserving homomorphisms. 
\end{prop}

\noindent Proof. Theorem \ref{ring18}  shows that the map $ \sigma \rightarrow \alpha _ \sigma, $ for $ \sigma \in D _ P, $ is 
onto the group $ k^* $. To show that the map is injective, let $ \sigma _1 $ and $ \sigma_2 $ be elements of $ D_P $ with $ \sigma_1 \ne  \sigma_2 $, and choose a point $ Q $ such that $ \sigma_1 Q \ne  \sigma_2 Q $. From the proof of theorem \ref{ring18}, we have $ \tau_{PQ}^{\alpha_{\sigma_i}} P = \sigma_i Q $. This shows that $ \alpha_ {\sigma_1} \ne \alpha_{\sigma_2} $. This argument is easily reversed; thus the  map is a bijection. The verification of the algebraic properties of an  isomorphism is straightforward. $\Box$

\begin{thm}
\label{pap4}
The following are equivalent.

(a) Postulate P.

(b) The division ring $ k $ of trace-preserving homomorphisms is commutative. 
\end{thm}

\noindent Proof. Let postulate P hold, let $ P $ be any point, and let $ \sigma _1 $ and $ \sigma _2 $ be any dilatations in the subgroup $ D _P $. Choose any two distinct lines $ l $ and $ m $ passing through $ P $; it follows that $ l \nparallel m $.  Choose any points $ Q $ and $ R $, on $ l $ and $ m $, each distinct from $ P $. Set $ Q' \equiv \sigma _1Q, $ set $ Q'' \equiv \sigma _2 \sigma _1 Q, $ set $ R' \equiv \sigma _2R, $ set $ R'' \equiv \sigma _1 \sigma _2 R, $ and set $ S \equiv \sigma _1 \sigma _2 Q. $ Since the dilatations are injective, these points are also distinct from $ P $. 
It follows from lemma \ref{dila23}  that $ Q', ~Q'',  $ and $ S $  lie on $ l $, while $  R'  $ and $ R'' $ lie on $ m $.  Thus $ Q + R' \parallel \sigma _1 Q + \sigma _1 R' = Q' + R'', $ and $ Q' + R \parallel \sigma _2 Q' + \sigma _2 R = Q'' + R'. $ 

Applying postulate P, we have $ Q'' + R'' \parallel Q + R \parallel \sigma _1 \sigma _2 Q + \sigma _1 \sigma _2 R = S + R'' $. It follows that both of the points $ Q'' $ and $ S $ lie on $ Q'' + R'', $ which is nonparallel to $ l $. Since the points  $ Q'' $ and $ S $ also  lie on $ l $, these points  are equal. Thus the dilatations $ \sigma _1 \sigma _2 $ and $ \sigma _2 \sigma _1 $ agree at the point $ Q $, and it follows from 
theorem \ref{dila14} that they  are equal. This shows that the subgroup $ D _P $ is commutative, and thus 
the group $ k^* $ is also commutative.

Conversely, let $ k^* $ be commutative, and let a Pappus configuration be given. Using axiom K2, construct dilatations $ \sigma _1 $ and $ \sigma _2 $ in the subgroup $ D _P $ such that $ \sigma _1Q = Q', $ and $ \sigma _2 R = R'. $ Now the above argument may be reversed. $\Box$

\section{Geometry based on a field}
\label{sec:geom} 

 Beginning now \label{geom1} with a given field $ k ,$ we  construct a geometry that satisfies all the axioms. The field $ k $ must possess certain special properties; these  are all constructively valid for the  real field $ \mathbb{R} $.\footnote{See [B, BB; Chapter 2].}

\begin{defn}
\label{geom2}
\textnormal{A \textit{Heyting field}\footnote{See [MRR; \S\,II.2].} is a field $ k $ with an inequality relation that is a tight apartness. Thus $ k $ satisfies the following conditions, where $ x $ and $ y $ are any elements of $ k $.\\
 \hspace*{5mm}(a) $ \neg (x \ne x ) $. \\
\hspace*{5mm}(b) If $  x \ne y,$ then 
$ y \ne x $.\\
 \hspace*{5mm}(c) If $ x \ne y ,$ then for any element $ z $, either $ z \ne x $ or $ z \ne y $.\\
 \hspace*{5mm}(d) If $ \neg (x \ne y) ,$ 
then $ x = y $.}  
\end{defn}

\begin{prop}
\label{geom3}
A Heyting field  $ k $ has the following  properties, 
where $ x $ and $ y $ are any elements of $ k $.

(e) $ x \ne 0 $ if and only if  $ x $ is a unit in $ k $.

(f) If $ xy \ne 0 ,$ then $ x \ne 0 $ and $ y \ne 0 $.
\end{prop}

\begin{prop}
\label{geom4}
The division ring of trace-preserving homomorphisms constructed in section \ref{sec:ring} satisfies all the conditions   for a Heyting field, and  the additional conditions (e) and (f), except for commutativity. 
\end{prop}

\noindent Proof. This follows from theorem \ref{ring20}. $ \Box $ \\

 \noindent \textit{Problem.} The classical theory [A] constructs a geometry based on an arbitrary division ring. Construct a geometry based on a division ring having all the conditions for a Heyting field, and the additional conditions (e) and (f), except for commutativity. 

\begin{defn}
\label{geom5} 
\textnormal{Let $ k $ be a Heyting field, and set $ \mathscr{P} _ k \equiv k^2 $. Let  $ P = (x,y) $ and 
$ Q = (z,w) $ be ``points'' in $ \mathscr{P} _ k $. We will write $ P = Q  $   if   $  x = z  $ and $  y = w ,$ and  $  P \ne Q  $   if either $   x \ne z  $ or $    y \ne w $. We will identify points and vectors in $ k ^2 ,$ and use vector notation; we will not use the notation $ l = P + Q $ in this section. For origin, 
set $ O \equiv (0,0) $.\\
\hspace*{5mm}For any point $ P $ in $ \mathscr{P} _ k $ and any vector $ A \ne 0 ,$ we define a set of points, which will be called a ``line,'' by
\begin{displaymath}
\overline {P,A} \equiv \{P + tA : t \in k\}
\end{displaymath}
The set of all such lines will be denoted $ \mathscr{L} _ k $. }
\end{defn} 

\begin{prop}
\label{geom6}
Definition \ref{geom5} defines a geometry 
$ \mathscr{G} _ k \equiv  ( \mathscr{P} _ k, \mathscr{L} _ k ) $ according to 
definition \ref{ax2}.
\end{prop}

 \noindent \textit{Note.} Theorem \ref{outside2}, a fundamental and quite essential result, ``If $ \neg (P \notin l) $, then $ P \in l$,'' was obtained in section \ref{sec:ax} only after introduction of the axioms; for $ \mathscr{G} _ k $ it may be proved directly by using coordinates. Let $ P = (x,y) $ be a point, and let $ l = \overline{R,A} $ be a line such that  $ \neg (P \notin l) $, where $ A = ( a_1, a_2 ) $. It suffices to consider the case in which $ a _1 \ne 0; $ thus, by a change in parameter, we may assume that $ R $ has coordinates $ ( x, r_2 ) $. Suppose that $ P \ne R, $ and let $ Q = R + tA $ be any point on $ l $. Either $ Q \ne P, $ or $ Q \ne R. $ In the second case, it follows from condition \ref{geom3}(f) that $ t \ne 0, $ and thus $ q _1 \ne x; $ this shows that $ Q \ne P $. Hence $ P \notin l, $ a contradiction; it follows that  $ P = R. $ Thus $ P \in l.$

\begin{lm}
\label{geom7}
Let $ l $ be a line.

(a) If $ l = \overline {P,A} $ and $ Q \in l ,$ then $ l $ may be written as $ l = \overline {Q,A} $.

(b) If $ l = \overline {P,A} $ and $ B = cA $ for some element $ c \ne 0 $ in $ k ,$ then $ l $ may be written as $ l = \overline {P,B} $.

(c) If $ P $ and $ Q $ are distinct points on $ l ,$ then $ l = \overline {P,Q-P} $. 
\end{lm}

\begin{lm}
\label{geom8}
Let $ l = \overline {P,A} $ and $ m = \overline {Q,B} $ be lines, with $ A = (a_1, a_2) $ and $ B = (b_1, b_2) $.

(a) $ l \nparallel m $ if and only if  $ a_1 b_2 \ne a_2 b_1 $.

(b) $ l \parallel m $ if and only if  $ a_1 b_2 = a_2 b_1 $.

(c) $ l \parallel m $ if and only if  $ B = cA $ for some element $ c \ne 0 $ in $ k $.
\end{lm}

\noindent Proof. Let $ l \nparallel m $; these lines then have a common point, and by a change of parameters we may assume that $ Q = P $. Since $ l \ne m ,$ we may assume that there exists a point $ R \in l $ with 
$ R \notin m $. Thus there exists an element $ t_1 \ne 0 $ in $ k $ such that
$ R = P + t_1A \ne P + tB ,$ for all 
$ t \in k $. We may assume that
$ b_1 \ne 0 $; take 
$ t \equiv t_1a_1b_1^ {-1} $.
Then $ A \ne a_1b_1^ {-1}B = (a_1, a_1b_1^{-1} b_2) ,$
and it follows that 
$ a_2 \ne a_1b_1^{-1} b_2 $. 
Thus* 
$ a_1b_2 \ne a_2b_1 $.\footnote{*In this section, portions of the proofs that utilize the commutativity of the field $ k $ will be marked with an asterisk.}

Now let $ a_1b_2 \ne a_2b_1 $. This condition  allows us to solve* the system resulting from the equation
$ P + tA = Q + uB $; thus $ l \cap m \ne \emptyset $. Choose a point $ R $ in $ l \cap m $; thus $ l = \overline{R,A} $ and $ m = \overline{R,B} $. Set
$ S \equiv R + A = (r_1 + a_1, r_2 + a_2) $; thus
$ S \in l $. Let $ T \equiv R + tB = (r_1 + tb_1, r_2 + tb_2) $ be any point on $ m $. Using condition \ref{geom2}(c), we have either $ tb_1b_2 \ne a_1b_2 $ or $ tb_1b_2 \ne a_2b_1 $. Thus, either  $ a_1 \ne tb_1 ,$ or* $ a_2 \ne tb_2 $; it follows that $ S \ne T $. Thus $ S \notin m ,$ and this shows that $ l \ne m $. Hence $ l \nparallel m $. $\Box$

\begin{thm}
\label{geom9}
The geometry $ \mathscr{G} _k $ satisfies the axioms in group $ \textbf{G} $. 
\end{thm}

\begin{thm}
\label{geom12}
The geometry $ \mathscr{G} _k $ satisfies axiom L2. 
\end{thm}

\noindent Proof. Let $ l = \overline{P,A} $ and 
$ m = \overline{P,B} $ be nonparallel lines; thus $ a_1b_2 \ne a_2b_1 $. Let $ n = \overline{Q,C} $ be any line; we may assume that $ c_2 \ne 0 $. Thus $ a_1b_2c_2 \ne a_2b_1c_2 $. Either $ c_1 a_2 b_2 \ne a_1 b_2 c_2 $ or $ c_1 a_2 b_2 \ne a_2 b_1 c_2 $. In the first case,* $ c_1a_2 \ne c_2a_1 ,$ and thus $ n \nparallel l $. In the second case,* $ c_1 b_2 \ne c_2 b_1 ,$ and thus $ n \nparallel m $. $\Box$

\begin{thm}
\label{geom13}
The geometry $ \mathscr{G} _k $ satisfies axiom L1.
\end{thm}

\noindent Proof. Let $ l \equiv \overline{P, A} $ and 
$ m \equiv \overline{P, B}$ be nonparallel lines, and let $ Q $ be a point distinct from $ P $. Set $ C \equiv Q - P $ and $ n \equiv \overline{P, C} $. Using axiom L2, we may assume that $ n \nparallel l $; thus $ a_1c_2 \ne a_2c_1 $. Now let $ R \equiv  P + tA $ be any point on $ l $. Either $ ta_1a_2 \ne a_1c_2 ,$ or $ ta_1a_2 \ne a_2c_1 $. Thus* either $ c_2 \ne ta_2 $ or $ c_1 \ne ta_1 $; hence $ Q \ne R $. This shows that $ Q \notin l $. $\Box$\\

\noindent \textit{Problem.} The\label{geom10} axioms in group \textbf{L} have been verified out of order, using axiom L2 to prove axiom L1. Prove\label{geom11} axiom L1 for $ \mathscr{G} _k $ directly; then use theorem \ref{ax36} to derive axiom L2.

\begin{thm}
\label{geom14}
The dilatations of the geometry  $\; \mathscr{G} _k $ are the maps $ \sigma $ defined by 
\begin{equation} 
\label{eq:dil} \sigma X \equiv eX + C  \qquad for \; all \; points  \; X 
\end{equation} 
where $ e \ne 0 $ is an element of $ k ,$ and $ C $ is any vector.

The translations of  $\; \mathscr{G} _k $ are the maps 
$ \tau _C ,$ where $ C $ is any vector, defined by
\begin{equation} 
\label{eq:tran} \tau _C X \equiv X + C  \qquad for \; all \; points  \; X 
\end{equation} 
The translation $ \tau _C $ is $ \ne 1 $ if and only if $ C \ne O $; in this case, the trace family of $ \tau _C $ is the pencil containing the line $ \overline{O,C} $. 
\end{thm}

\noindent Proof. Let $ \sigma $ be a map defined by 
(\ref{eq:dil}); it is clearly injective. Let $ P  $ and $ Q $ be distinct points, let $ l $ be the line through $ P $  and $ Q ,$ and let $ l' $ the line through $ \sigma P $ and $ \sigma Q $. Then $ l' = \overline{\sigma P,\sigma Q-\sigma P} = \overline{\sigma P,e(Q-P)} = \overline{\sigma P,Q-P} \parallel \overline{P,Q-P} = l
$. This shows that $ \sigma $ is a dilatation.

Now let $ \sigma $ be any  dilatation. Set $ C \equiv \sigma O ,$ and choose any point $ U $ distinct from $ O $. Since $ \sigma $ takes the points $ O $ and $ U $ into the points $ C $ and $ \sigma U ,$ the lines $ \overline{O,U} $ and $ \overline{C, \sigma U - C} $ are parallel. Thus there exists an element $ e \ne 0 $ in $ k $ such that $ \sigma U - C = eU $. Let $ \sigma' $ be the dilatation defined by (\ref{eq:dil}), using these values for $ e $ and $ C $. Since $ \sigma $ and $ \sigma' $ agree at the points $ O $ and $ U ,$ it follows from theorem \ref{dila14} that $ \sigma = \sigma' $. This shows that the maps defined in (\ref{eq:dil}) include all dilatations.

The traces, if any, of the dilatation $ \tau  _C $ have the form $ \overline{X,C} $. Since these lines are all parallel, $ \tau  _C $ is a translation. Now let $ \tau $ be an arbitrary translation, choose any point $ P ,$ and set $ C \equiv \tau P - P $. Since $ \tau _C $ agrees with $ \tau $ at the point $ P ,$  it follows from theorem \ref{tran17} that $ \tau = \tau _C $. Thus the maps defined in (\ref{eq:tran}) include all translations. $\Box$

\begin{thm}
\label{geom15}
The geometry $ \mathscr{G} _k $ satisfies the axioms in group \textbf{K}. 
\end{thm}

\noindent Proof. Given points $ P $ and $ Q ,$ set 
$ C \equiv Q - P $;  the translation $ \tau _C $  satisfies axiom K1. Now let $ Q $ and $ R $ be points collinear with, and distinct from, the origin $ O $. There exists an element $ e \ne 0 $ in $ k $ such that 
$ R = eQ $; define a dilatation by $ \sigma X \equiv  eX $. Thus proposition \ref{coor3A} applies, and axiom K2 is valid. $\Box$\\

Since the geometry $ \mathscr{G} _k $ now satisfies the axioms in all three groups, it is a constructive Desarguesian plane.

\begin{thm}
\label{geom16}
The trace-preserving homomorphisms of the geometry 
$ \mathscr{G}_k $ are the 
maps $ \alpha _x ,$ for all $ x \in k ,$ defined by

\begin{equation} \label{eq:tph} 
\tau_C ^{\alpha_x} \equiv \tau_{xC} \qquad for \; all \; translations \; \tau _C 
\end{equation} 
Furthermore, $ \alpha _x  \ne 0 $ if and only 
if $ x \ne 0 $.
\end{thm}

\noindent Proof. Let $ x \in k $. The algebraic condition for a trace-preserving homomorphism is easily verified for $ \alpha _x   $. Now let 
$ \tau _C $ be any translation, and let $ l $ be a trace of $ \tau _C ^ { \alpha _x } $. It follows that $ \tau _{xC} \ne 1 $; thus $ x \ne 0 $ and $  C \ne O $. Then 
$ l_1 \equiv \overline{O,xC} $ is a trace of 
$ \tau _C ^ { \alpha _x } ,$ and 
$ l_2 \equiv \overline{O,C} $ is a trace of $ \tau _C $. Since $ l \parallel l_1 $ 
and $ l_1 \parallel l_2 ,$ it follows that 
$ l \in t (\tau _C) $. This shows 
that $ \alpha _x $ is a trace-preserving homomorphism.

Now let $ \alpha $ be any trace-preserving homomorphism. Set 
$ C = (1,0) ,$ and $ D \equiv (x,y) \equiv \tau _C ^ \alpha O $. Suppose that $ y \ne 0 $. Then $ \overline{O,D} $ is a trace of $ \tau _C ^ \alpha ,$ and since $ \overline{O,C} $ is a trace of $ \tau _C ,$  these lines are parallel. It follows that   
$ D = dC $ for some  $ d \in k,$ and 
$ y = 0 $, a contradiction. Hence $ y = 0 $. Thus $ \tau _C ^ \alpha O = (x,0) = \tau _{xC}O = \tau _C ^ { \alpha _x} O $.
Since the translations $ \tau _C ^ \alpha $ and $ \tau _C ^ { \alpha _x} $ agree at the origin $ O ,$ it follows from 
theorem \ref{tran17} that they are equal. Since $ \alpha $ and $ \alpha _x $ agree at the translation $ \tau _C \ne 1 ,$ it follows from corollary \ref{ring14}(c) that they are equal. This shows that the maps defined in (\ref{eq:tph}) include all trace-preserving homomorphisms. $\Box$

\begin{thm}
\label{geom18}
The set $  \overline{k} $ of trace-preserving homomorphisms of the translation group $ T $ in the geometry $ \mathscr{G} _k $ is a field, isomorphic with the given field $ k $ under the map $ x \rightarrow \alpha _x $ of $ k $ onto $  \overline{k} $. Let $ \mathscr{G}_k $ be coordinatized by the field $ \overline{k} $ as in section \ref{sec:coor}, using the point $ (0,0) $ as origin in theorem \ref{coor6}, with 
$ \tau _1 \equiv \tau _{(1,0)} $ and $ \tau _2 \equiv \tau _{(0,1)} $.

If, using the field 
$ \overline{k} ,$ a point 
$ P = (x,y) $ in $ \mathscr{P} _k $ is assigned the coordinates $ ( \xi , \eta  ) ,$ then 
$ \xi = \alpha _x $ and 
$ \eta = \alpha  _y $. Thus the two coordinate systems correspond under the isomorphism 
$ x \rightarrow \alpha _x $ between the two fields.
\end{thm}

\noindent Proof. The algebraic properties for an isomorphism are easily verified. In 
theorem \ref{coor6}, the 
coordinates  $ ( \xi , \eta ) $ are assigned to the point $ P = (x,y) $ according to the rule 
$ \tau_{OP} = 
\tau_1 ^\xi \tau_2 ^\eta,$ while  
$ \tau_{OP} $ is in the present section denoted 
$ \tau_P $. Now, $ \tau _P = \tau _{(x,y)} = 
\tau _{(x,0)} \tau _{(0,y)} = 
\tau _1 ^{\alpha _x} \tau _2 ^{\alpha _y} $. Thus it follows from the uniqueness shown in theorem \ref{coor5} that 
$ \xi = \alpha _x $ and 
$ \eta = \alpha  _y $. $ \Box $

\section{The real plane}
\label{sec:real} 
  
 The constructive properties of the real field $ \mathbb{R}$ ensure that it is a Heyting field;\footnote{The various properties are found in [B, BB; Chapter 2].} thus  it follows from section  \ref{sec:geom} that $ \mathbb{R}^{2} $  is a Desarguesian plane. 

The field $ \mathbb{R} $ has additional structures compared to an arbitrary field, especially order and metric. This raises the  possibility of other choices for the principal relations on the plane $ \mathbb{R} ^{2}$. Clearly, the primitive relation $ P \ne Q $  is equivalent to the condition $ \rho (P,Q) > 0 $. Theorem \ref{real1} will show that the principal relation, $ P \notin l $, as  given in definition \ref{ax-2b}, is also equivalent to conditions involving the additional structures.  

A subset $ F $ of a constructive metric space $ (M, \rho ) $ is \textit{located in M} if the distance $ \rho (x,F) \equiv \mathrm{inf}_{y \in F} \, \rho (x,y) $  may be determined for any point $ x $ in $ M $. Any line $ l $ on the real plane $ \mathbb{R}^{2} $ is a located subset.\footnote{This follows from the results in [B, BB; Chapter 4].} An equation may be  found for any line, as noted in the comment following the proof of theorem \ref{coor6}. 

\begin{thm}  
\label{real1}
Let $ \rho $ be the usual metric 
on $ \mathbb{R}^{2} $. Let 
$ P = (x_0, y_0) $ be any point, and let $ l $ be a line with equation $ ax + by + c = 0. $ Then the following are equivalent:

(a) $ P \notin l $

(b) $ \rho(P,l) > 0 $

(c) $ ax_0 + by_0 + c \ne 0 $
\end{thm}

\noindent Proof. If $ \rho (P,l) > 0,$ then for any point $ Q $ on $ l $, we have $ \rho (P,Q) \geq \rho (P,l) > 0 $, and thus $ P \ne Q. $ This shows that $ P \notin l $. Conversely, if $ P \notin l, $ then $ P \ne Q $ for all points $ Q $ on $ l $, and it follows 
from [B; Chapter 6, Lemma 7]\footnote {For an alternative proof of this lemma, see [M1, Lemma 5.4].} that $ \rho (P,l) > 0.$ 

 Now let $ ax_0 + by_0 + c \ne 0. $ Let $ Q = (x_1,y_1) $ be any point on $ l $; thus $  ax_1 + by_1 + c = 0. $ It follows that $  0 < | a (x_0 - x_1) + b (y_0 - y_1) | \leq | a | | x_0 - x_1 | + | b | | y_0 - y_1 |. $ Thus at least one of the last two terms is positive, and $ P \ne Q. $ This shows that $ P \notin l. $

Finally, let $ P \notin l. $ Let $ m $ be the line $ y = y_0 $, and let $ n $  be the line  $ x = x_0 $. It follows from lemma \ref{geom8} that $ m \nparallel n $.
Using axiom L2, we may assume that $ l \nparallel m; $ it then follows that $ a \ne 0. $ Set $ Q \equiv l \cap m ; $ thus  $ Q $ has coordinates of the form $ Q  = (x_1,y_0)$. Since $ P \notin l, $ we have 
$ P \ne Q, $ and thus $ x_0 \ne x_1. $ Since $ Q \in l $, we have $ ax_1 + by_0 + c = 0, $ and thus 
$ ax_0 + by_0 + c = a (x_0 - x_1) \ne 0. $ $\Box$\\

 \noindent \textit{Problem.} For an arbitrary Heyting field $ k $, extend the  part of theorem \ref{real1} involving conditions (a) and (c) to the geometry $ \mathscr{G} _ k $ constructed in section \ref{sec:geom}.\\

\section{Brouwerian counterexamples}
\label{sec:brou} 

To determine the specific nonconstructivities in the classical theory, and the points at which modification is required, we use \textit{Brouwerian counterexamples}, in conjunction with \textit{omniscience principles.} A Brouwerian counterexample contains a proof that a given statement implies an omniscience principle. In turn, an omniscience principle, taken with full constructive meaning, would provide solutions, or significant information, for a large number of well-known unsolved problems.\footnote{This method was introduced by L. E. J. Brouwer in 1908 to demonstrate that the content of mathematics is placed in jeopardy by use of the \textit{principle of the excluded middle}. For a discussion of Brouwer's critique and the reaction of the mathematical community, see [S, page 319ff]. For more information concerning Brouwerian counterexamples, and other omniscience principles, see [BR], [MRR], [M1], and [M2].} 

For example, the results of an effort to find a counterexample to the Goldbach conjecture might be recorded as a binary sequence: set $ a_n = 0 $ if you verify the conjecture up through $  n, $  and set $ a_n = 1 $ when you find a counterexample $\leq n$. Given an arbitrary binary sequence $ (a_n) $, the \textit{limited principle of omniscience (LPO)}\footnote{LPO and LLPO were introduced by Brouwer, and given the current names by Errett Bishop.} provides either a proof that $ a_n = 0 $ for all $ n $, or a finite routine for constructing an integer $ n $ with $ a_n = 1 $; this would settle the Goldbach problem, along with many other unsolved problems. No one has, nor is it conceivable that anyone will ever find, such a general principle. While humans may discover proofs that settle certain individual questions, only an omniscient being would claim to possess a finite routine for predicting the outcome of an arbitrary infinite search. Thus the principle LPO is considered nonconstructive. 

Although the omniscience principles are usually stated in terms of binary sequences, these sequences may be used to construct corresponding real numbers; this results in the following equivalent formulations for the principal omniscience principles:

\textit{Limited principle of omniscience (LPO).} For any real number $ c \geq 0, $ either $ c = 0 $ or $ c > 0. $

 \textit{Weak limited principle of omniscience (WLPO).} For any real number $ c \geq  0, $ either $ c =0 $ or $ \neg (c = 0). $

 \textit{Lesser limited principle of omniscience (LLPO).} For any real number $ c, $ either 
$ c \leq 0 $ or $ c \geq 0. $

 \textit{Limited principle of existence (LPE).}\footnote{The principle LPE is usually called \textit{Markov's principle (MP).} Although accepted in the Markov school of recursive function theory, this  principle is nonconstructive according to the strict constructivism introduced by Bishop. No strictly constructive algorithm validating this principle is known, and it is unlikely that such an algorithm will ever be found. Markov's principle asserts a general finite routine: Given an infinite binary sequence, and a proof that it is contradictory that each term is 0, MP constructs a positive integer $ n $ such that the $n^{\mathrm{th}}$ term is 1. For more information concerning Markov's principle, see [BR].} For any real number ~$ c \geq 0, $~ if $ \neg (c=0) $, then $ c > 0. $ 

A statement will be considered \textit{nonconstructive} \label{noncx} if it implies one of these omniscience principles. The examples in this section all take place on the  real plane $ \mathbb{R}^{2}. $ 

\begin{ex}
\label{brouJ}\label{brouM}
\textnormal{The following statements are nonconstructive.\\
 \hspace*{5mm}(i) \textit{Given any points $ P $ and $ Q $, either $ P = Q $ or $ P \ne Q $}.\\
 \hspace*{5mm}(ii) \textit{Given any point $ P $ and any line l, either $ P \in l $ or $ P \notin l $}.\\
 \hspace*{5mm}(iii) \textit{Given any lines $ l $ and $ m $, either $ l \parallel m $ or $ l \nparallel m $.}}
 \end{ex}
 
 \noindent Let $ c \geq 0 $ be a real number.  Set $ P \equiv (0,c), $  set $ Q \equiv (0,0), $ let $ l $ be the line $ y = 0,$ and let $ m $ be the line $ y = cx $.  Each statement implies LPO.

\begin{ex}
\label{brouA}
\textnormal{The following statements are nonconstructive.\\
 \hspace*{5mm}(i) \textit{If $ \neg (P=Q), $ then $ P \ne Q $}.\\
 \hspace*{5mm}(ii) \textit{If $ \neg (P \in l), $ then $ P \notin l $}.}
 \end{ex}
 
 \noindent Let $ c \geq 0 $ be a real number such that $ \neg(c = 0). $ Set $ P \equiv (0,c), $ set $ Q = (0,0), $ and let $ l $ be the line $ y = 0. $ Each statement implies LPE.

\begin{ex}
\label{brouB}
\textnormal{The following statement is nonconstructive.\\
 \hspace*{5mm}\textit{If the lines $ l $ and $ m $ are parallel, then either $ l = m $ or $ l \cap m = \emptyset $.}\footnote{This is the converse of proposition \ref{ax32}.}}
\end{ex}

\noindent Let $ c \geq 0 $ be a real number, and define the lines $ l $ and $ m $ by $ y = 0 $ and $ y = c $. Suppose that $ l \nparallel m. $ Then these lines have a common point, so $ c = 0 $, and they are distinct, so $ c > 0 $. Hence $ l \parallel m. $ The  statement implies WLPO. 

\begin{ex}
\label{brouC}\label{brouD}
\textnormal{The following statements are nonconstructive.\\
 \hspace*{5mm}(i) \textit{If $ l $ and $ m $ are lines with $ \neg (l \parallel m), $ then $ l \nparallel m $}.\\
 \hspace*{5mm}(ii) \textit{If the lines $ l $ and $ m $ have a unique point in common, then $ l \nparallel m $.}}
\end{ex}
 
\noindent Let $ c \geq 0 $ be a real number such that 
$ \neg (c = 0) $ and let $ l $ and $ m $ be the lines $ y = 0 $ and $ y = cx $.\\
\hspace*{5mm}(i) Suppose that $ l \parallel m; $ then $ l = m, $ and thus $ c = 0, $ a contradiction. Hence $ \neg (l \parallel m). $ The statement implies that  $ l \nparallel m, $ and thus $ l \ne m. $ It follows that one of the lines contains a point that is outside the other line. In one case, we have a point $ (x, cx) \notin l; $ thus $ (x, cx) \ne (x, 0). $ In the other case, we have a point $ (x, 0) \notin m; $ thus $ (x, 0) \ne (x, cx). $ In either case, $ cx \ne 0 $, and thus
 $ c \ne 0. $  This shows that statement (i) implies LPE.\\
\hspace*{5mm}(ii) Let $ P = (x,y) $ be any common point; thus $ cx = 0. $ Suppose that $ P \ne (0,0); $ then $ x \ne 0, $ and thus $ c = 0, $ a contradiction. Hence $ P = (0,0). $ Thus $ l $ and $ m $ have a unique common point. Statement (ii) implies that $ l \nparallel m $; it has been shown in part (i) that this implies LPE.

\begin{ex}
\label{brouE}
\textnormal{The following statements are nonconstructive.\\
 \hspace*{5mm}(i) \textit{Any given dilatation is either the identity, or distinct from the identity.}\\
\hspace*{5mm}(ii) \textit{Any given dilatation  either has a fixed point, or has no fixed point.}\\
\hspace*{5mm}(iii) \textit{Any given translation is either the identity, or has no fixed point.}} 
\end{ex}

\noindent Let $ c\geq 0 $ be a real number. Using theorem \ref{geom14}, define a  dilatation 
$ \sigma $ by $ \sigma X \equiv X + (c,0). $ The first statement implies LPO; the last two statements each imply WLPO. 

\begin{ex}
\label{brouL}
\textnormal{Assume for the moment that the definition of ``dilatation'' were to allow the ``degenerate'' case, in which all points map onto a single point. Then the following statement  [A, Theorem 2.3] is nonconstructive.\footnote{This example relates to a comment following definition \ref{dila2}.}\\
 \hspace*{5mm}\textit{Any given dilatation is either degenerate or (weakly) injective.}}  
\end{ex}

\noindent Let $ c \geq 0$ be a real number, and consider the map $ X \rightarrow cX $. The statement  
implies (WLPO) LPO.

\begin{ex}
\label{brouF}
\textnormal{The following statement [A, Theorem 2.12] is nonconstructive.\footnote{This example relates to the comment preceding corollary \ref{ring14}.} \\
\hspace*{5mm}\textit{Let $ \tau $ be a translation, and let $ \alpha $ be a trace-preserving homomorphism. If $ \tau ^ \alpha = 1, $ then either $ \alpha = 0 $ or $ \tau = 1 $.}}
\end{ex} 

\noindent Let $ c $ be any real number. Set $ d \equiv \textnormal{max} \{c,0\}, $ and $ e \equiv \textnormal{min} \{c,0\};$ thus $ de=0. $  Set $ \tau \equiv \tau _{(e,0)}, $ and set $ \alpha \equiv \alpha _d. $ Thus $ \tau ^ \alpha = \tau _ {(de,0)} = \tau _O = 1 $. If $ \alpha = 0 $, then $ d = 0 $, and $ c \leq 0 $. If $ \tau = 1 $, then $ e = 0 $ , and  $ c \geq 0 $.  Thus the statement implies LLPO.

\begin{ex}
\label{brouK}
\textnormal{The following statement is nonconstructive.\footnote{This example relates to lemma \ref{tran12}.}\\
\hspace*{5mm}\textit{For any translation $ \tau, $ there exists a pencil of lines $ \pi $ such 
that the trace family of \,$ \tau $ is contained in 
$  \pi. $}}
\end{ex}

\noindent Let $ c, ~d, $ and $ e $ be as in example  \ref{brouF}, and let $ \tau $ be the translation defined by $ \tau X \equiv X + (d, e). $ Use the statement to choose a line $ l $ such that $ t ( \tau ) \subseteq \pi _l. $ It follows from axiom L2 that $ l $ is either nonparallel to the line $ y = 0 $ or nonparallel to the line $ x = 0. $ In the first case, suppose that $ c > 0; $ then $ d > 0, ~e = 0, $ and the line $ y = 0 $ is a trace of $ \tau, $ a contradiction. Hence $ c \leq 0. $ The second case is similar. Thus the statement implies LLPO.

\begin{ex} 
\label{brouH}\label{brouI}
\textnormal{Weakening the definition of ``nonparallel'' adopted in definition \ref{ax8} is not feasible. It must follow from the definition that nonparallel lines have a common point, and axiom L1 must be allowed. The following statements are nonconstructive.\\
\hspace*{5mm}(i) \textit{If $ l \ne m $ and $ \neg (l\cap m = \emptyset), $ then $ l \nparallel m $}.\\
\hspace*{5mm}(ii) \textit{If $ \neg (l = m) $ and $ l \cap m \ne \emptyset , $ then $ l \nparallel m $}.\\ 
\hspace*{5mm}(iii) \textit{If ~$ \neg (l \cap m \ne \emptyset $ implies $ l = m), $ 
then $ l \nparallel m $}.\footnote{This example relates to the comment following proposition \ref{ax34}.}} 
\end{ex} 

\noindent Let $ c \geq 0 $ be a real number such that $ \neg (c = 0), $ and let $ l $ be the line $ y = 0. $

(i) Let $ m $ be the line $ y = cx + 1; $ it is clear that $ l \ne m. $ Now suppose that $ l \cap m = \emptyset . $ Suppose further that $ c > 0; $ then $ l \cap m \ne \emptyset, $ a contradiction. Hence $ c = 0, $ a contradiction. Thus $ \neg (l  \cap m = \emptyset ).$ The statement implies that $ l \nparallel m; $ thus the lines have a common point $ (x, y), $ and it follows that $ c > 0. $ Thus statement (i) implies LPE.

(ii) Let $ m $ be the line $ y = cx. $ The statement implies that $ l \nparallel m $.  It follows from axiom L1 that $ (1,0) \notin m, $ and hence $ c > 0. $ Thus statement (ii) implies LPE.

(iii) Let $ m $ be the line in part (i). Suppose that the implication holds. Suppose further that $ c > 0; $ then $ l \cap m \ne \emptyset, $ and it follows that $ l = m, $ a contradiction. Hence $ c = 0, $  a contradiction. Thus the implication is contradictory. Statement (iii) implies that $ l \nparallel m; $ it has been shown in part (i) that this implies LPE.\\

\normalsize
\setlength{\parindent}{0pt}
\section*{References} 
\label{sec:ref}

[A] E. Artin, \textit{Geometric algebra.} Interscience Publishers, Inc., New York-London, 1957.\\	

[B] E. Bishop, \textit{Foundations of constructive analysis.} McGraw-Hill Book Co., New York-Toronto-London, 1967.\\ 

[BB] E. Bishop and D. Bridges, \textit{Constructive analysis.} Springer-Verlag, Berlin, 1985.\\

[BM] D. Bridges and R. Mines, What is constructive mathematics? Math. Intel. 6 (1984), 32-38.\\

[BR] D. Bridges and F. Richman, \textit{Varieties of constructive mathematics.} Cambridge University Press, Cambridge, 1987.\\

[BV] D. Bridges and L. V\^\i\c{t}\u{a}, \textit{Techniques of constructive analysis.} Springer, New York, 2006.\\

[D1] D. van Dalen, Extension problems in intuitionistic plane projective geometry I, II, Indag. Math. 25 (1963), 349-383.\\

[D2] D. van Dalen, ``Outside'' as a primitive notion in constructive projective geometry, Geom. Dedicata 60 (1996), 107-111.\\

[H] D. Hilbert, \textit{Grundlagen der Geometrie.} Verlag B. G. Teubner, Leipzig, 1899.\\

[He1] A. Heyting, Zur intiuitionistischen Axiomatik der projektiven Geometrie, Math. Ann. 98 (1927), 491-538.\\

[He2] A. Heyting, Axioms for intuitionistic plane affine geometry, in L. Henkin, P. Suppes, A. Tarski (eds), \textit{The axiomatic method, with special reference to geometry and physics: Proceedings of an international symposium held at the University of California, Berkeley, December 26, 1957 - January 4, 1958}. North-Holland, Amsterdam, 1959, 160-173.\\

[L] D. Li, Using the prover ANDP to simplify orthogonality, Ann. Pure Appl.
Logic 124 (2003), 49-70.\\

[LJL] D. Li, P. Jia, X. Li, Simplifying von Plato's axiomatization of constructive apartness geometry, Ann. Pure Appl. Logic 102 (2000), 1-26.\\

[LV] M. Lombard, R. Vesley, A common axiom set for classical and intuitionistic plane geometry, Ann. Pure Appl. Logic 95 (1998), 229-255.\\

[M1] M. Mandelkern, \textit{Constructive continuity.} Memoirs Amer. Math. Soc. 42 (1983), nr. 277.\\

[M2] ---------------, Limited omniscience and the Bolzano-Weierstrass principle, Bull. London Math. Soc. 20 (1988), 319-320.\\

[MRR] R. Mines, F. Richman, and W. Ruitenburg, \textit{A course in constructive algebra.} Springer-Verlag, New York, 1988.\\ 

[P1] J. van Plato, The axioms of constructive geometry, Ann. Pure Appl. Logic 76 (1995), 169-200.\\

[P2] J. van Plato, A constructive theory of ordered affine geometry, Indag.
Math. 9 (1998), 549-562.\\

[R] F. Richman, Existence proofs, Amer. Math.  Monthly 106 (1999), 303-308.\\ 

[S] G. Stolzenberg, Review of E. Bishop, \textit{Foundations of constructive analysis.} Bull. Amer. Math. Soc. 76 (1970), 301-323.\\

\hspace*{5mm} \\
\hspace*{5mm} \\

Department of Mathematics \\
New Mexico State University \\\\
mandelkern@zianet.com \\
www.zianet.com/mandelkern \\
Postal address: \\
Mark Mandelkern \\
5259 Singer Road \\
Las Cruces NM 88007-5566 \\ 
USA

\end{document}